\documentclass{amsart}
\usepackage{times,amssymb,amscd,graphics,axodraw,verbatim}
\usepackage[all]{xy}
\usepackage[stdtext,vcentermath]{myyoungtab}

\numberwithin{equation}{section}

\newtheorem{prop}{Proposition}
\newtheorem{theorem}[prop]{Theorem}
\newtheorem{corollary}[prop]{Corollary}
\newtheorem{lemma}[prop]{Lemma}
\newtheorem{conjecture}[prop]{Conjecture}

\theoremstyle{definition}
\newtheorem{definition}[prop]{Definition}
\newtheorem{example}[prop]{Example}
\newtheorem{remark}[prop]{Remark}

\numberwithin{prop}{section}

\newcommand{\A}{\mathcal{A}}
\newcommand{\bin}[2]{\left( \begin{array}{c} #1\\ #2 \end{array} \right)}
\newcommand{\cc}{\mathrm{cc}}
\newcommand{\Complex}{\mathbb{C}}
\newcommand{\Conf}{\mathrm{C}}
\newcommand{\Confb}{\overline{\Conf}}
\newcommand{\CST}{\mathrm{CST}}
\newcommand{\SCST}{\mathrm{SCST}}
\newcommand{\Dt}{D}
\newcommand{\ev}{\mathrm{ev}}
\newcommand{\HH}{\mathcal{H}}
\newcommand{\Hil}{\mathcal{H}}
\newcommand{\Ih}{\widehat{I}}
\newcommand{\inner}[2]{\langle #1\,,\,#2\rangle}
\newcommand{\Jb}{\overline{J}}
\newcommand{\la}{\lambda}
\newcommand{\La}{\Lambda}
\newcommand{\lb}{\mathrm{lb}}
\newcommand{\lh}{\mathrm{lh}}
\newcommand{\lm}{\la^-}
\newcommand{\ls}{\mathrm{ls}}
\newcommand{\lt}{\widetilde{\ell}}
\newcommand{\Mb}{\overline{M}}
\newcommand{\mt}{m^t}

\newcommand{\nub}{\overline{\nu}}
\newcommand{\Path}{\mathcal{P}}
\newcommand{\Pathb}{\overline{\Path}}
\newcommand{\Phib}{\overline{\Phi}}
\newcommand{\pr}{\mathrm{pr}}
\newcommand{\R}{\mathbb{R}}
\newcommand{\RC}{\mathrm{RC}}
\newcommand{\RCb}{\overline{\RC}}

\newcommand{\rcls}{i}
\newcommand{\rclb}{j}
\newcommand{\rev}{\mathrm{rev}}
\newcommand{\rk}{\mathrm{rk}}
\newcommand{\row}{\mathrm{row}}
\newcommand{\SA}{\mathcal{SA}}
\newcommand{\sbf}{\vec{s}}
\newcommand{\sigbf}{\vec{\sigma}}
\newcommand{\trP}{\mathrm{tr}_{\Path}}
\newcommand{\trRC}{\mathrm{tr}_{\RC}}
\newcommand{\qbin}[2]{\genfrac{[}{]}{0pt}{}{#1}{#2}}
\newcommand{\qbins}[2]{{\textstyle\genfrac{[}{]}{0pt}{}{#1}{#2}}}
\newcommand{\wt}{\mathrm{wt}}
\newcommand{\x}{x}
\newcommand{\Xb}{\overline{X}}
\newcommand{\Z}{\mathbb{Z}}

\begin{document}

\title[$X=M$ Theorem]{$X=M$ Theorem:\\
{\small Fermionic formulas and rigged configurations under review}}

\author[A.~Schilling]{Anne Schilling}
\address{Department of Mathematics, University of California, One Shields
Avenue, Davis, CA 95616-8633, U.S.A.}
\email{anne@math.ucdavis.edu}
\urladdr{http://www.math.ucdavis.edu/\~{}anne}
\thanks{\textit{Date:} October 2005}
\thanks{Partially supported by the NSF grants DMS-0200774 and DMS-0501101.}

\begin{abstract}
We give a review of the current status of the $X=M$ conjecture. Here $X$ stands for
the one-dimensional configuration sum and $M$ for the corresponding fermionic formula.
There are three main versions of this conjecture: the unrestricted, the classically 
restricted and the level-restricted version. We discuss all three versions and
illustrate the methods of proof with many examples for type $A_{n-1}^{(1)}$.
In particular, the combinatorial approach via crystal bases and rigged configurations
is discussed. Each section ends with a conglomeration of open problems.
\end{abstract}

\maketitle

\tableofcontents

\newpage

\section{Introduction}
The \textit{$X=M$ Conjecture} asserts the equality between the generating function
of highest weight tensor product crystal elements graded by the energy function
and the fermionic formula~\cite{HKOTT:2001,HKOTY:1999,O:2005}. This article concerns the 
\textit{$X=M$ Theorem}, or more precisely, those cases in which the $X=M$ Conjecture 
has been proven. We describe the method of proof which uses the combinatorics of 
crystal bases and rigged configurations.
We mostly focus on type $A_{n-1}^{(1)}$, but many of the constructions have 
analogues for other affine Kac--Moody algebras $\mathfrak{g}$. Instead of providing
all details of the proofs, we illustrate the main concepts via examples. Each
section ends with a conglomeration of open problems.

The \textit{fermionic formula} is a $q$-analogue of the tensor product multiplicity
$[\bigotimes_j W_{s_j}^{(r_j)}, V_\lambda]$, where $W_s^{(r)}$ is a $U_q(\mathfrak{g})$
Kirillov--Reshetikhin module indexed by a Dynkin node $r$ and $s\in \Z_{>0}$,
and $V_\lambda$ is the irreducible highest weight $U_q(\overline{\mathfrak{g}})$-module 
with highest weight $\lambda$. Here $\overline{\mathfrak{g}}$ is the finite-dimensional
classical algebra inside the affine Kac--Moody algebra $\mathfrak{g}$. Alternatively, 
since the procedure of taking the crystal limit does not change tensor product 
multiplicities, we can view the fermionic formula as a $q$-analogue of 
$[\bigotimes_j B^{r_j,s_j},B(\lambda)]$, where $B^{r,s}$ is the Kirillov--Reshetikhin 
crystal and $B(\lambda)$ is the finite-dimensional
highest weight crystal indexed by the dominant weight $\lambda$. Instead of labeling
the fermionic formula by $B=\bigotimes_j B^{r_j,s_j}$ and $B(\lambda)$, we use
the multiplicity array $L=(L_s^{(r)})$ and $\lambda$, where $L_s^{(r)}$ denotes the
number of tensor factors $B^{r,s}$ in $B$. For type $A_{n-1}^{(1)}$ the fermionic
formula is then given by
\begin{equation}\label{eq:fermi intro}
  \Mb(L,\la;q)=\sum_{\nu\in\Confb(L,\la)} q^{\cc(\nu)}
   \prod_{(a,i)\in\HH} \qbin{p_i^{(a)}+m_i^{(a)}}{m_i^{(a)}}.
\end{equation}
Here $\Confb(L,\lambda)$ is the set of admissible $(L,\lambda)$-configurations,
$\HH=I\times \Z_{>0}$ with $I=\{1,2,\ldots,n-1\}$, $m_i^{(a)}$ is the particle number 
and $p_i^{(a)}$ is the vacancy number. The precise definition of the various quantities 
is given in section~\ref{sec:RC}. The $q$-binomial coefficient is defined as
\begin{equation*}
\qbin{p+m}{m}=\frac{(q)_{p+m}}{(q)_p(q)_m}
\end{equation*}
for $p,m\in\Z_{\ge 0}$ and zero otherwise, where $(q)_m=(1-q)(1-q^2)\cdots (1-q^m)$.

The $q$-binomial coefficient $\qbin{p+m}{m}$ is the generating function of 
partitions in a box of size $p\times m$. Using this interpretation, 
equation~\eqref{eq:fermi intro} can be rewritten in solely combinatorial terms as
\begin{equation*}
\Mb(L,\lambda;q)=\sum_{(\nu,J)\in \RCb(L,\lambda)} q^{\cc(\nu,J)},
\end{equation*}
where $\RCb(L,\lambda)$ is the set of \textit{rigged configurations} as defined in 
section~\ref{sec:RC}. 
The one-dimensional
configuration sum $\Xb(B,\lambda;q)$ is the generating function of highest weight paths 
$\Pathb(B,\lambda)$ of weight $\lambda$ weighted by the energy function $\Dt$
\begin{equation*}
\Xb(B,\lambda;q)=\sum_{b\in \Pathb(B,\lambda)} q^{\Dt(b)}.
\end{equation*}
The $X=M$ conjecture~\cite{HKOTT:2001,HKOTY:1999} asserts that
\begin{equation}\label{eq:X=M}
\Xb(B,\lambda;q)=\Mb(L,\lambda;q)
\end{equation}
for all affine Kac--Moody algebras $\mathfrak{g}$.

The $X=M$ conjecture can be proved by establishing a statistics preserving bijection 
$\Phib:\Pathb(B,\la)\to \RCb(L,\la)$ between the set of paths and the set of rigged 
configurations. More precisely, $\Phib$ should have the property that 
$\Dt(b)=\cc(\Phib(b))$ for all $b\in \Pathb(B,\lambda)$. For 
$B=\bigotimes_j B^{1,\mu_j}$ of type $A_{n-1}^{(1)}$ such a bijection was given by 
Kerov, Kirillov and Reshetikhin~\cite{KKR:1986,KR:1988}. In fact,
in this case the set of paths $\Pathb(B,\lambda)$ is in bijection with the set
of semi-standard Young tableaux $\mathrm{SSYT}(\lambda,\mu)$ of shape $\lambda$ and
content $\mu=(\mu_1,\mu_2,\ldots)$, and the energy function corresponds to the cocharge 
of Lascoux and Sch\"utzenberger~\cite{LS:1978}. The bijection of Kerov, Kirillov and 
Reshetikhin~\cite{KKR:1986,KR:1988} is a bijection between semi-standard Young tableaux
and rigged configurations and yields a fermionic formula for the Kostka--Foulkes
polynomials. In~\cite{KSS:2002}, this bijection was generalized to
$B=\bigotimes_j B^{r_j,s_j}$ of type $A_{n-1}^{(1)}$. In this case the set of paths
$\Pathb(B,\lambda)$ is in bijection with Littlewood--Richardson tableaux and the
bijection was in fact formulated as a bijection between Littlewood--Richardson
tableaux and rigged configurations. For other types such bijections have also been
given in special cases. In summary to date the following cases have been proven:
\begin{itemize}
\item $B=\bigotimes_j B^{r_j,s_j}$ of type $A_{n-1}^{(1)}$~\cite{KSS:2002};
\item $B=\bigotimes_j B^{1,s_j}$ of all nonexceptional types~\cite{OSS:2003,SS:2004};
\item $B=\bigotimes_j B^{r_j,1}$ for type $D_n^{(1)}$~\cite{S:2005}. 
\end{itemize}
An important technique in studying fermionic formulas of nonsimply-laced types are
virtual crystals and virtual rigged configuration~\cite{OSS:2003a,OSS:2003b}.

In this paper we provide a review of the bijective approach to the $X=M$ conjecture.
We will mostly restrict our attention to type $A_{n-1}^{(1)}$ and set up the bijection
between crystals and rigged configurations (rather than tableaux and rigged 
configurations).

The correspondence between the two combinatorial sets can be understood in terms
of two approaches to solvable lattice models and their associated spin chain systems: 
the \textit{Bethe Ansatz}~\cite{Bethe:1931} and the 
\textit{corner transfer matrix method}~\cite{Baxter:1982}.

In his 1931 paper~\cite{Bethe:1931}, Bethe solved the Heisenberg spin chain
based on the string hypothesis which asserts that the eigenvalues of
the Hamiltonian form certain strings in the complex plane as the size
of the system tends to infinity. The Bethe Ansatz has been applied to many
further models proving completeness of the Bethe vectors.
The eigenvalues and eigenvectors of the Hamiltonian are 
indexed by rigged configuration. However, numerical studies indicate 
that the string hypothesis is not always true~\cite{ADM:1992}.
The corner transfer matrix (CTM) method was introduced by Baxter and labels
the eigenvectors by one-dimensional lattice paths. It turns out that these
lattice paths have a natural interpretation in terms of Kashiwara's
crystal base theory~\cite{Kash:1990}, namely as highest weight crystal
elements in a tensor product of finite-dimensional crystals.

\begin{figure}
\resizebox{10cm}{!}{\includegraphics{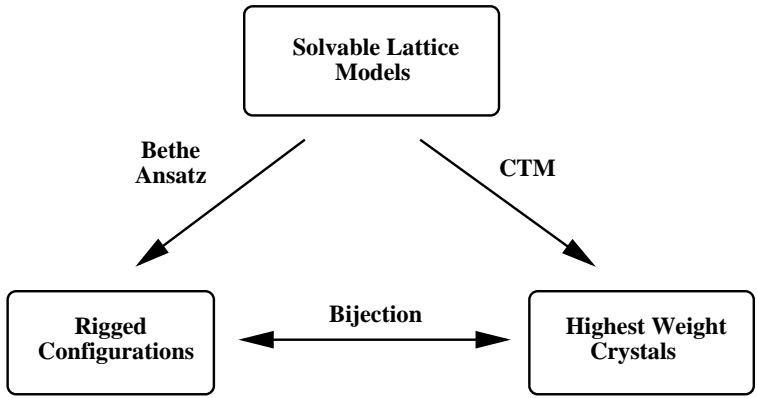}}
\caption{Bethe Ansatz versus corner transfer matrix method (CTM)\label{fig:bethe ctm}}
\end{figure}

Even though neither the Bethe Ansatz nor the corner transfer matrix
method are mathematically rigorous, they suggest that there should be
a bijection between the two index sets, namely rigged configurations on the
one hand and highest weight crystal elements on the other hand.
This is schematically indicated in Figure~\ref{fig:bethe ctm}.
As explained above, the generating function of rigged configurations leads fermionic 
formulas. Fermionic formulas can be interpreted as explicit expressions for the 
partition function of the underlying physical models which reflect the particle 
structure. For more details regarding the physical background of fermionic formulas
see~\cite{KKMM:1993a,KKMM:1993b,HKOTT:2001}.

The $X=M$ conjecture can be generalized in two different ways: to the 
\textit{level-restricted} and the \textit{unrestricted} case. Both of these cases 
will also be reviewed in this paper in the case of type $A_{n-1}^{(1)}$.

The set of paths $\Pathb(B,\lambda)$ is defined as the set of all $b\in B$ of
weight $\lambda$ that are highest weight with respect to the classical crystal
operators. The Kirillov--Reshetikhin crystals are affine crystals and have the
additional crystal operators $e_0$ and $f_0$, which can be used to define 
level-restricted paths. Hence it is natural to consider the generating functions of
level-restricted paths, giving rise to a level-restricted version of $X$.
The corresponding set of level-restricted rigged configurations was considered
in~\cite{SS:2001}. The notion of level-restriction is also very important in the
context of restricted-solid-on-solid (RSOS) models in statistical
mechanics~\cite{Baxter:1982} and fusion models in conformal field 
theory~\cite{V:1988}. The one-dimensional configuration sums of
RSOS models are generating functions of level-restricted paths
(see for example~\cite{ABF:1984,DJKMO:1987,JMO:1987}). The structure constants of the 
fusion algebras of Wess--Zumino--Witten conformal field theories
are exactly the level-restricted analogues of the tensor product multiplicities
$X(B,\lambda;1)$ or Littlewood--Richardson coefficients as shown by 
Kac~\cite[Exercise 13.35]{Kac:1990} and Walton~\cite{Wal:1990a,Wal:1990b}. $q$-Analogues 
of these level-restricted Littlewood--Richardson coefficients in terms of ribbon 
tableaux were proposed in ref.~\cite{FLOT:1998}.

Rigged configurations corresponding to highest weight crystal paths are only the 
tip of an iceberg. In~\cite{S:2005a} the definition of rigged configurations was extended
to all crystal elements in types $ADE$ by the explicit construction of a crystal 
structure on the set of \textit{unrestricted} rigged configurations. 
The equivalence of the crystal structures on rigged configurations and crystal paths 
together with the correspondence for highest weight vectors yields the equality 
of generating functions in analogy to~\eqref{eq:X=M}. Denote the unrestricted
set of paths and rigged configurations by $\Path(B,\la)$ and $\RC(L,\la)$, respectively.
The corresponding generating functions are unrestricted one-dimensional configuration 
sums or $q$-supernomial coefficients. A direct bijection 
$\Phi:\Path(B,\la)\to \RC(L,\la)$ for type $A_{n-1}^{(1)}$ along the lines 
of~\cite{KSS:2002} is constructed in~\cite{DS:2004,DS:2005}.

The paper is organized as follows. In section~\ref{sec:bethe} we present the Bethe
Ansatz for the spin $1/2$ XXX Heisenberg chain which first gave rise to rigged configurations.
In section~\ref{sec:X} we review the one-dimensional configuration sums and set the notation 
used in this article. The corresponding fermionic formulas
for the classically restricted, unrestricted and level-restricted cases are
subject of sections~\ref{sec:classical}, \ref{sec:unrestricted} and~\ref{sec:level},
respectively. In particular for the $\Xb=\Mb$ case, we introduce rigged configurations 
and fermionic formulas in section~\ref{sec:RC}, define certain splitting operations on 
crystals and rigged configurations in sections~\ref{sec:operations crystals} 
and~\ref{sec:operations rc}, which are necessary for the bijection $\Phib$ between
paths and rigged configurations of section~\ref{sec:bij}. Section~\ref{sec:prop} 
features many of the properties of $\Phib$. For the unrestricted version of 
the $X=M$ theorem, we define the crystal structure on rigged configurations
in section~\ref{sec:rc crystal}. A characterization of unrestricted rigged
configurations is given in section~\ref{sec:lower bound} which is used in 
section~\ref{sec:fermi} to derive the fermionic formula. The affine crystal
operators on rigged configurations are given in section~\ref{sec:affine}.
Section~\ref{sec:level} deals with the level-restricted version of the $X^\ell=M^\ell$
theorem. Level-restricted rigged configurations are introduced in section~\ref{sec:level rc}
and the corresponding fermionic formula is derived in section~\ref{sec:level fermi}.
Each section ends with some open problems.

\subsection*{Acknowledgments}
I would like to thank Atsuo Kuniba and Masato Okado for organizing the workshop
``Combinatorial Aspect of Integrable Systems'' at the Research Institute for Mathematical 
Sciences in Kyoto in July 2004 for which this review was written.

\section{Bethe Ansatz and rigged configurations}
\label{sec:bethe}

In this section we discuss the algebraic Bethe Ansatz
for the example of the spin $1/2$ XXX Heisenberg chain and show
how rigged configurations arise. Further details can be found
in~\cite{Fadeev:1998,S:2003}.

The spin $1/2$ XXX \textbf{Heisenberg chain} is a one-dimensional quantum 
spin chain on $N$ sites
with periodic boundary conditions. It is defined on the Hilbert space
$\Hil_N = \bigotimes_{n=1}^N h_n$ where in this case $h_n=\Complex^2$
for all $n$. Associated to each site is a local spin variable
$\sbf = \frac{1}{2} \sigbf$ where
\begin{equation*}
\sigbf=(\sigma^1,\sigma^2,\sigma^3)=
\left( \left(\begin{array}{cc} 0&1\\ 1&0\end{array}\right),
\left( \begin{array}{cc} 0&-i\\ i&0 \end{array} \right),
\left( \begin{array}{cc} 1&0\\ 0&-1 \end{array} \right) \right)
\end{equation*}
are the Pauli matrices.
The spin variable acting on the $n$-th site is given by
\begin{equation*}
\sbf_n = I \otimes \cdots \otimes I \otimes \sbf \otimes I
 \otimes \cdots \otimes I
\end{equation*}
where $I$ is the identity operator and $\sbf$ is in the $n$-th
tensor factor. We impose periodic boundary conditions
$\sbf_n = \sbf_{n+N}$.

The \textbf{Hamiltonian} of the spin $1/2$ XXX model is
\begin{equation*}
H_N = J \sum_{n=1}^N \left(\sbf_n \cdot \sbf_{n+1} - \frac{1}{4}\right).
\end{equation*}
Our goal is to determine the eigenvectors and eigenvalues of 
$H_N$ in the antiferromagnetic regime $J>0$ in the limit when
$N\to \infty$.

The main tool is the \textbf{Lax operator} $L_{n,a}(\la)$, also called the local
transition matrix. It acts on $h_n\otimes \Complex^2$ where $\Complex^2$ is an
auxiliary space and is defined as
\begin{equation*}
L_{n,a}(\la) = \la I_n \otimes I_a + i \sbf_n \otimes \sigbf_a.
\end{equation*}
Here $I_n$ and $I_a$ are unit operators acting on $h_n$ and the
auxiliary space $\Complex^2$, respectively; $\la$ is a complex parameter,
called the spectral parameter. Writing the action on the auxiliary space as
a $2\times 2$ matrix, we have
\begin{equation}\label{eq:L aux}
L_{n}(\la) = \left( \begin{array}{cc} \la + i s_n^3 & i s_n^-\\
i s_n^+ & \la - i s_n^3 \end{array} \right)
\end{equation}
where $s_n^\pm = s_n^1 \pm i s_n^2$.

The crucial fact is that the Lax operator satisfies commutation relations
in the auxiliary space $V=\Complex^2$. Altogether there are 16 relations
which can be written compactly in tensor notation.
Given two Lax operators $L_{n,a_1}(\la)$ and $L_{n,a_2}(\mu)$ defined
in the same quantum space $h_n$, but different auxiliary spaces
$V_1$ and $V_2$, the products $L_{n,a_1}(\la)L_{n,a_2}(\mu)$ and
$L_{n,a_2}(\mu)L_{n,a_1}(\la)$ are defined on the triple tensor
product $h_n\otimes V_1\otimes V_2$.
There exists an operator $R_{a_1,a_2}(\la-\mu)$ defined on
$V_1\otimes V_2$ such that
\begin{equation}\label{eq:com rel}
R_{a_1,a_2}(\la-\mu) L_{n,a_1}(\la)L_{n,a_2}(\mu)
= L_{n,a_2}(\mu)L_{n,a_1}(\la) R_{a_1,a_2}(\la-\mu).
\end{equation}
Explicitly, the $R$-matrix $R_{a_1,a_2}(\la)$ is given by
\begin{equation*}
R_{a_1,a_2}(\la) = \left( \la +\frac{i}{2}\right) I_{a_1}\otimes I_{a_2}
 + \frac{i}{2} \sigbf_{a_1} \otimes \sigbf_{a_2}.
\end{equation*}

Geometrically, the Lax operator $L_{n,a}(\la)$ can be interpreted 
as the transport between sites $n$ and $n+1$ of the quantum
spin chain. Hence
\begin{equation*}
T_{N,a}(\la) = L_{N,a}(\la) \cdots L_{1,a}(\la)
\end{equation*}
is the monodromy around the circle (recall that we assume periodic
boundary conditions). In the auxiliary space write
\begin{equation*}
T_{N}(\la) = \left( \begin{array}{cc}
A(\la) & B(\la)\\
C(\la) & D(\la)
\end{array} \right)
\end{equation*}
with entries in the full Hilbert space $\Hil_N$. From (\ref{eq:com rel})
it is clear that the \textbf{monodromy matrix} satisfies the
following commutation relation
\begin{equation}\label{eq:com rel mon}
R_{a_1,a_2}(\la-\mu)T_{N,a_1}(\la)T_{N,a_2}(\mu)
= T_{N,a_2}(\mu)T_{N,a_1}(\la)R_{a_1,a_2}(\la-\mu).
\end{equation}

Let $\omega_n=\left( \begin{array}{c} 1\\0\end{array} \right)$. 
In the auxiliary space the Lax operator is triangular on $\omega_n$
\begin{equation*}
L_n(\la)\omega_n = \left( \begin{array}{cc} \la+\frac{i}{2} & *\\
0 & \la-\frac{i}{2} \end{array} \right) \omega_n
\end{equation*}
where $*$ stands for an for us irrelevant quantity. This follows
directly from (\ref{eq:L aux}). On the Hilbert
space $\Hil_N$ we define $\Omega= \bigotimes_n \omega_n$ so that
\begin{equation*}
T_N(\la) \Omega = \left( \begin{array}{cc} \alpha^N(\la) & * \\
0 & \delta^N(\la) \end{array} \right) \Omega
\end{equation*}
where $\alpha(\la)=\la+\frac{i}{2}$ and $\delta(\la)=\la-\frac{i}{2}$.
Equivalently this means that
\begin{eqnarray*}
C(\la) \Omega &=& 0\\
A(\la) \Omega &=& \alpha^N(\la) \Omega\\
D(\la) \Omega &=& \delta^N(\la) \Omega
\end{eqnarray*}
so that $\Omega$ is an eigenstate of $A(\la)$ and $D(\la)$
and hence also of $t_N(\la)=A(\la)+D(\la)$.

The claim is that the other eigenvectors of $t_N(\la)$ are of the form
\begin{equation*}
\Phi(\la,\La) = B(\la_1)\cdots B(\la_n) \Omega.
\end{equation*}
The lambdas $\La=\{\la_1,\ldots,\la_n\}$ 
satisfy a set of algebraic relations, called the
\textbf{Bethe equations}, which can be derived from~\eqref{eq:com rel mon}
\begin{equation}\label{eq:bethe 1}
\left(\frac{\la+\frac{i}{2}}{\la-\frac{i}{2}}\right)^N
=\prod_{\stackrel{\la'\in \La}{\la'\neq \la}} 
 \frac{\la-\la'+i}{\la-\la'-i}
\end{equation}
where $\la\in\La=\{\la_1,\ldots,\la_n\}$.

Suggested by numerical analysis, it is assumed that in the limit
$N\to\infty$ the $\la$'s form strings.
This hypothesis is called the \textbf{string hypothesis}.
A string of length $\ell=2M+1$, where $M$ is an integer or half-integer
depending on the parity of $\ell$, is a set of $\la$'s of the form
\begin{equation*}
\la^M_{jm} = \la^M_j + im
\end{equation*}
where $\la^M_j\in\R$ and $-M\le m\le M$ is integer or half-integer
depending on $M$.
The index $j$ satisfies $1\le j\le m_\ell$ where $m_\ell$
is the number of strings of length $\ell$.
A decomposition of $\{\la_1,\ldots,\la_n\}$ into strings is called
a \textbf{configuration}. Each configuration is parametrized by $\{m_\ell\}$.
It follows that
\begin{equation*}
\sum_\ell \ell m_\ell = n.
\end{equation*}

Now take \eqref{eq:bethe 1} and multiply over a string
\begin{eqnarray}\label{eq:inter}
\prod_{m=-M}^M &&
 \left( \frac{\la^M_j+i(m+\frac{1}{2})}{\la^M_j+i(m-\frac{1}{2})}
 \right)^N\nonumber\\
&=&\prod_{m=-M}^M \prod_{\stackrel{M',j',m'}{(M',j',m')\neq (M,j,m)}}
\frac{\la^M_j-\la^{M'}_{j'}+i(m-m'+1)}{\la^M_j-\la^{M'}_{j'}+i(m-m'-1)}.
\end{eqnarray}
Many of the terms on the left and right cancel so that this equation
can be rewritten as
\begin{equation}\label{eq:bethe 2}
e^{i N p_M(\la^M_j)} = \prod_{\stackrel{M',j'}{(M',j')\neq (M,j)}} 
 e^{i S_{MM'}(\la_j^M-\la_{j'}^{M'})},
\end{equation}
in terms of the momentum and scattering matrix
\begin{eqnarray*}
e^{ip_M(\la)} &=& \frac{\la+i(M+\frac{1}{2})}{\la-i(M+\frac{1}{2})}\\
e^{i S_{MM'}(\la)} &=& \prod_{m=|M-M'|}^{M+M'}
\frac{\la+im}{\la-im}\cdot \frac{\la+i(m+1)}{\la-i(m+1)}.
\end{eqnarray*}
Taking the logarithm of \eqref{eq:bethe 2} using the branch cut
\begin{equation*}
\frac{1}{i}\ln \frac{\la+ia}{\la-ia} = \pi -2\arctan \frac{\la}{a}
\end{equation*}
we obtain
\begin{equation}\label{eq:Q}
2N \arctan \frac{\la^M_j}{M+\frac{1}{2}} = 2\pi Q^M_j
 + \sum_{\stackrel{M',j'}{(M',j')\neq (M,j)}}
 \Phi_{MM'}(\la^M_j-\la^{M'}_{j'}),
\end{equation}
where
\begin{equation*}
\Phi_{MM'}(\la) = 2 \sum_{m=|M-M'|}^{M+M'}
 \left(\arctan \frac{\la}{m}+\arctan\frac{\la}{m+1}\right).
\end{equation*}
The first term on the right is absent for $m=0$.
Here $Q^M_j$ is an integer or
half-integer depending on the configuration.

In addition to the string hypothesis, we assume that the $Q^M_j$
classify the $\la$'s uniquely: $\la^M_j$ increases if $Q^M_j$
increases and in a given string no $Q^M_j$ coincide. As we will see
shortly with this assumption one obtains the correct number
of solutions to the Bethe equations \eqref{eq:bethe 1}.

Using $\arctan \pm \infty = \pm \frac{\pi}{2}$ we obtain from
(\ref{eq:Q}) putting $\la^M_j=\infty$
\begin{equation*}
Q^M_\infty = \frac{N}{2} -\bigl(2M+\frac{1}{2}\bigr)\bigl(m_{2M+1}-1\bigr)
- \sum_{M'\neq M}\bigl(2\min(M,M')+1\bigr) m_{2M'+1}.
\end{equation*}
Since there are $2M+1$ strings in a given string of length $2M+1$,
the maximal admissible $Q^M_{\mathrm{max}}$ is
\begin{equation*}
Q^M_{\mathrm{max}} = Q^M_\infty -(2M+1)
\end{equation*}
where we assume that if $Q^M_j$ is bigger than $Q^M_{\mathrm{max}}$
then at least one root in the string is infinite and hence all
are infinite which would imply $Q^M_j=Q^M_\infty$.

With the already mentioned assumption that each admissible set of
quantum number $Q^M_j$ corresponds uniquely to a solution of the
Bethe equations we may now count the number of Bethe vectors.
Since $\arctan$ is an odd function and by the assumption about the
monotonicity we have
\begin{equation*}
-Q^M_{\mathrm{max}}\le Q^M_1<\cdots <Q^M_{m_{2M+1}}\le Q^M_{\mathrm{max}}.
\end{equation*}
Hence defining $p_\ell$ as
\begin{equation*}
p_\ell=N-2\sum_{\ell'} \min(\ell,\ell')m_{\ell'}
\end{equation*}
so that
\begin{equation*}
p_\ell+m_\ell=2Q^M_{\max}+1 \qquad \mbox{with $\ell=2M+1$}.
\end{equation*}
With this the number of Bethe vectors with configuration
$\{m_\ell\}$ is given by
\begin{equation*}
Z(N,n|\{m_\ell\}) = \prod_{\ell\ge 1} 
 \bin{p_\ell+m_\ell}{m_\ell}
\end{equation*}
where $\bin{p+m}{m}=(p+m)!/p!m!$ is the binomial
coefficient. The total number of Bethe vectors is
\begin{equation}\label{eq:count}
Z(N,n) = \sum_{\stackrel{\{m_\ell\}}{\sum_\ell \ell m_\ell=n}}
 \prod_{\ell\ge 1} \bin{p_\ell+m_\ell}{m_\ell}.
\end{equation}

It should be emphasized that the derivation of (\ref{eq:count})
given here is not mathematically rigorous. Besides the various
assumptions that were made we also did not worry about possible
singularities of (\ref{eq:inter}). However, (\ref{eq:count}) indeed 
yields the correct number of Bethe vectors.

To interpret (\ref{eq:count}) combinatorially let us view
the set $\{m_\ell\}$ as a partition $\nu$. A partition is a
set of numbers $\nu=(\nu_1,\nu_2,\ldots)$ such that $\nu_i\ge \nu_{i+1}$
and only finitely many $\nu_i$ are nonzero. The partition has
part $i$ if $\nu_k=i$ for some $k$. The size of partition $\nu$
is $|\nu|:=\nu_1+\nu_2+\cdots$. In the correspondence
between $\{m_\ell\}$ and $\nu$, $m_\ell$ specifies the number of parts
of size $\ell$ in $\nu$. For example, if $m_1=1$, $m_2=3$, $m_4=1$ and
all other $m_\ell=0$ then $\nu=(4,2,2,2,1)$.

It is well-known (see e.g. \cite{Andrews:1976}) that 
$\bin{p+m}{m}$ is the number of partitions
in a box of size $p\times m$, meaning, that the partition cannot
have more than $m$ parts and no part exceeds $p$. Let
$\RCb(N,n)$ be the set of all \textbf{rigged configurations} $(\nu,J)$ defined 
as follows. $\nu$ is a partition of size $|\nu|=n$ and $J$ is a
set of partition where $J_\ell$ is a partition in a box of
size $p_\ell\times m_\ell$. Then (\ref{eq:count}) can be rewritten as 
\begin{equation}\label{eq:Z rc}
Z(N,n) = \sum_{(\nu,J)\in \RCb(N,n)} 1.
\end{equation}

\begin{example}\label{ex:rig}
Let $N=5$ and $n=2$. Then the following is the set of
rigged configuration $\RCb(5,2)$
\begin{eqnarray*}
& \yngrc(2,{1,1}) & \qquad \yngrc(1,{1,1},1,{1,1})\\
& \yngrc(2,{0,1}) & \qquad \yngrc(1,{1,1},1,{0,1})\\
& \mbox{} & \qquad \yngrc(1,{0,1},1,{0,1}).
\end{eqnarray*}
The underlying partition on the left is (2) and on the right (1,1).
The partitions $J_\ell$ attached to part length $\ell$ is specified
by the first number next to each part. For example, the partition $J_1$ for the top
rigged configuration on the right is (1,1) whereas for the one in the
middle and bottom is $J_1=(1)$ and $J_1=\emptyset$, respectively.
The numbers to the right of part $\ell$ is $p_\ell$.
\end{example}
The rigged configurations introduced in this section correspond to the
algebra $A_1$. In section~\ref{sec:RC}, we introduce rigged configurations
for the type $A_{n-1}$ algebras and also define a statistics $\cc$ which
turns~\eqref{eq:Z rc} into a polynomial in $q$.

\section{One-dimensional configuration sums and crystals}
\label{sec:X}
One-dimensional configuration sums are generating functions of crystal elements.
A detailed account on crystals can for example be found 
in~\cite{HKOTT:2001,HKOTY:1999,Kash:1990,O:2005}. 
Here we review the main definitions to fix our notation. We restrict ourselves
to crystals associated to $\mathfrak{g}$ of type $A_{n-1}^{(1)}$.

A crystal path is an element in the tensor product of crystals 
$B=B^{r_k,s_k}\otimes B^{r_{k-1},s_{k-1}}\otimes \cdots \otimes B^{r_1,s_1}$,
where $B^{r,s}$ is the Kirillov--Reshetikhin crystal labeled by
$r\in I=\{1,2,\ldots,n-1\}$ and $s\in \Z_{>0}$.
As a set the crystal $B^{r,s}$ of type $A_{n-1}^{(1)}$ is the set of all column-strict
Young tableaux of shape $(s^r)$ over the alphabet $\{1,2,\ldots,n\}$. 
Kashiwara~\cite{Kash:1990} introduced the notion of crystals and crystal graphs as a 
combinatorial means to study representations of quantum algebras. In particular, there 
are Kashiwara operators $e_i,f_i$ defined on the elements in $B^{r,s}$ for $0\le i<n$.

We first focus on $e_i,f_i$ when $i\in I$.
Let $b=b_k\otimes b_{k-1}\otimes \cdots \otimes
b_1\in B^{r_k,s_k}\otimes B^{r_{k-1},s_{k-1}}\otimes \cdots \otimes B^{r_1,s_1}$.
Let $\row(b)=\row(b_k)\row(b_{k-1})\ldots\row(b_1)$ be the concatenation of the
row reading words of $b$. For a fixed $i$, consider the subword of $\row(b)$
consisting of $i$'s and $(i+1)$'s only. Successively bracket all pairs $i+1\; i$.
What is left is a subword of the form $i^a (i+1)^b$. Define
\begin{equation*}
\begin{split}
e_i(i^a (i+1)^b)&=\begin{cases} i^{a+1} (i+1)^{b-1} & \text{if $b>0$}\\
0 & \text{otherwise} \end{cases}\\
f_i(i^a (i+1)^b)&=\begin{cases} i^{a-1} (i+1)^{b+1} & \text{if $a>0$}\\
0 & \text{otherwise.} \end{cases}
\end{split}
\end{equation*}

\begin{example}
Let $b=\young(12,23)\otimes \young(23) \otimes \young(1,3)$. Then $\row(b)=23122331$,
$e_2(\row(b))=23122231$ and $e_1(\row(b))=23112331$, so that
\begin{equation*}
\begin{split}
e_1(b)&=\young(11,23)\otimes \young(23) \otimes \young(1,3)\\
e_2(b)&=\young(12,23)\otimes \young(22) \otimes \young(1,3).
\end{split}
\end{equation*}
\end{example}

There are several sets of paths that will play an important role in the following.
For a composition of nonnegative integers $\la$, the set of 
\textbf{unrestricted paths} is defined as
\begin{equation*}
\Path(B,\la)=\{b\in B\mid \wt(b)=\la\}.
\end{equation*}
Here $\wt(b)=(w_1,\ldots,w_n)$ is the weight of $b$ where $w_i$ counts the number
of letters $i$ in $b$. For a partition $\la=(\la_1,\la_2,\ldots,\la_n)$, the set of
\textbf{classically restricted paths} is defined as
\begin{equation*}
\Pathb(B,\la)=\{b\in B\mid \wt(b)=\la, \quad \text{$e_i(b)=0$ for all $1\le i<n$}\}.
\end{equation*}
\begin{example}
For $B=B^{1,1}\otimes B^{2,2}\otimes B^{3,1}$ of type $A^{(1)}_3$ and 
$\la=(3,3,1,1)$ the path
\begin{equation*}
b\;=\; \young(2) \otimes \young(11,24) \otimes \young(1,2,3)
\end{equation*}
is in $\Pathb(B,\la)$.
\end{example}  

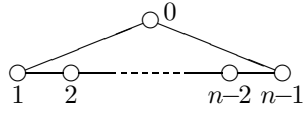
\begin{figure}
\begin{picture}(106,40)(-5,-5)
\multiput( 0,0)(20,0){2}{\circle{6}}
\multiput(80,0)(20,0){2}{\circle{6}} \put(50,20){\circle{6}}
\multiput( 3,0)(20,0){2}{\line(1,0){14}}
\multiput(63,0)(20,0){2}{\line(1,0){14}}
\multiput(39,0)(4,0){6}{\line(1,0){2}}
\put(2.78543,1.1142){\line(5,2){44.429}}
\put(52.78543,18.8858){\line(5,-2){44.429}}
\put(0,-5){\makebox(0,0)[t]{$1$}}
\put(20,-5){\makebox(0,0)[t]{$2$}}
\put(80,-5){\makebox(0,0)[t]{$n\!\! -\!\! 2$}}
\put(100,-5){\makebox(0,0)[t]{$n\!\! -\!\! 1$}}
\put(55,20){\makebox(0,0)[lb]{$0$}}
\end{picture}
\caption{Dynkin diagram for $A_{n-1}^{(1)}$ \label{fig:A}}
\end{figure}

There is a third set of level-restricted paths. The definition of these
paths requires the affine Kashiwara crystal operators $e_0$ and $f_0$. 
The affine Dynkin diagram of type $A_{n-1}^{(1)}$ has a circular symmetry,
which looks like a cycle with vertices labeled by $\Z/ n\Z$ (see Figure~\ref{fig:A}). 
The affine crystal $B^{r,s}$ also has such a symmetry, where the map 
$i\mapsto i+1 \pmod{n}$ on the vertices of the Dynkin diagram corresponds to the 
\textbf{promotion} operator $\pr$. Then the action of $e_0$ and $f_0$ is given by
\begin{equation}\label{eq:affine}
\begin{split}
e_0 &= \pr^{-1} \circ e_1 \circ \pr,\\
f_0 &= \pr^{-1} \circ f_1 \circ \pr.
\end{split}
\end{equation}
The promotion operator is a bijection $\pr:B\to B$ such that the following
diagram commutes for all $i\in \Ih=\{0,1,\ldots,n-1\}$
\begin{equation}\label{eq:pr com}
\begin{CD}
B@>{\pr}>> B \\
@V{f_i}VV @VV{f_{i+1}}V \\
B @>>{\pr}> B
\end{CD}
\end{equation}
and such that for every $b\in B$ the weight is rotated
\begin{equation}\label{eq:weight rotation}
\inner{h_{i+1}}{\wt(pr(b))}=\inner{h_i}{\wt(b)}.
\end{equation}
Here subscripts are taken modulo $n$. 

The promotion operator can be defined combinatorially using jeu de taquin~\cite{Sh:2002}.
Let $t\in B^{r,s}$ be a rectangular tableau of shape $(s^r)$. Delete all letters 
$n$ from $t$ and use jeu de taquin to slide the boxes into the empty spaces
until the shape of the new tableau is of skew shape $(s^r)/(\mu_n)$ where $\mu_n$
is the number of $n$ in $t$. Add one to all letters and fill the empty spaces 
by $1$s. The result is $\pr(t)$.
\begin{example}
Suppose $n=5$ and let
\begin{equation*}
t= \young(123,234,455).
\end{equation*}
Then removing the letters 5 and performing jeu de taquin, we obtain
\begin{equation*}
\young({}{}2,133,244).
\end{equation*}
Hence
\begin{equation*}
\pr(t)=\young(113,244,355).
\end{equation*}
\end{example}
\begin{example}
Take $t$ from the previous example. Then
\begin{equation*}
f_0(t)=\pr^{-1}\circ f_1 \circ \pr(t)= \young(112,233,445).
\end{equation*}
\end{example}
The set of \textbf{level-$\ell$ restricted paths} is now defined as
\begin{equation*}
\Path^{\ell}(B,\la)=\{b\in B\mid \wt(b)=\la, \quad \text{$e_i(b)=0$ for all $1\le i<n$},
\quad e_0^{\ell+1}(b)=0\}.
\end{equation*}

There exists a crystal isomorphism $R:B^{r,s}\otimes B^{r',s'} \to B^{r',s'} \otimes 
B^{r,s}$, called the \textbf{combinatorial $R$-matrix}. Combinatorially it is given 
as follows. Let $b\in B^{r,s}$ and $b'\in B^{r',s'}$. The product $b\cdot b'$ of two 
tableaux is defined as the Schensted insertion of $b'$ into
$b$. Then $R(b\otimes b')=\tilde{b}'\otimes \tilde{b}$ is the unique pair of tableaux
such that $b\cdot b'=\tilde{b}'\cdot\tilde{b}$.

The \textbf{local energy function} $H:B^{r,s}\otimes B^{r',s'}\to \Z$ is defined as
follows. For $b\otimes b'\in B^{r,s}\otimes B^{r',s'}$, $H(b\otimes b')$ is the number
of boxes of the shape of $b\cdot b'$ outside the shape obtained by concatenating $(s^r)$ 
and $({s'}^{r'})$.

\begin{example}
For 
\begin{equation*}
b\otimes b'= \young(12,24) \otimes \young(1,3,4)
\end{equation*}
we have 
\begin{equation*}
b\cdot b' = \young(113,224,4) = \young(1,2,4) \cdot \young(13,24) = \tilde{b}'\cdot
\tilde{b}.
\end{equation*}
so that
\begin{equation*}
R(b\otimes b')=\tilde{b}'\otimes\tilde{b}=\young(1,2,4) \otimes \young(13,24).
\end{equation*}
Since the concatentation of $\yng(2,2)$ and $\yng(1,1,1)$ is $\yng(3,3,1)$, the local
energy function $H(b\otimes b')=0$.
\end{example}

Now let $B=B^{r_k,s_k}\otimes \cdots\otimes B^{r_1,s_1}$ be a $k$-fold tensor
product of crystals. The \textbf{tail energy function} $\Dt:B\to \Z$ is given by
\begin{equation*}
  \Dt = \sum_{1\le i<j\le k} H_{j-1} R_{j-2} \dotsm R_{i+1} R_i,
\end{equation*}
where $H_i$ (resp. $R_i$) is the local energy function (resp. combinatorial $R$-matrix)
acting on the $i$-th and $(i+1)$-th tensor factors.

\begin{definition}\label{def:X}
The \textbf{one-dimensional configuration sum} is the generating function of 
the corresponding set of paths graded by the tail energy function
\begin{equation*}
\begin{split}
X(B,\la;q)&=\sum_{b\in \Path(B,\la)} q^{\Dt(b)},\\
\Xb(B,\la;q)&=\sum_{b\in \Pathb(B,\la)} q^{\Dt(b)},\\
X^{\ell}(B,\la;q)&=\sum_{b\in \Path^{\ell}(B,\la)} q^{\Dt(b)}.
\end{split}
\end{equation*}
The generating functions are called unrestricted, classically restricted and
level-restricted one-dimensional configuration sums or generalized Kostka
polynomials, respectively.
\end{definition}

\subsection{Open Problems}
\begin{itemize}
\item For types other than $A_{n-1}^{(1)}$, the existence of the Kirillov--Reshetikhin
crystals $B^{r,s}$ has been conjectured in~\cite{HKOTT:2001,HKOTY:1999}. The existence
of $B^{r,s}$, their combinatorial structure and properties are not yet well-understood
in general.
For the nonsimply-laced cases, the theory of virtual crystals~\cite{OSS:2003a,OSS:2003b}
can be employed to obtain the combinatorial structure of these crystal in terms
of the simply-laced cases.
\item For types other than $A_{n-1}^{(1)}$, a combinatorial construction of $R$ and $\Dt$
needs to be given.
\end{itemize}

\section{$\Xb=\Mb$}\label{sec:classical}

In this section we consider the $\Xb=\Mb$ theorem for type $A_{n-1}^{(1)}$, which 
was proven in~\cite{KSS:2002}. We begin by defining the fermionic formula $\Mb(L,\la;q)$
in section~\ref{sec:RC} and then describe the bijection $\Phib:\Pathb(B,\la)\to \RC(L,\la)$
and its properties in sections~\ref{sec:bij} and~\ref{sec:prop}.

\subsection{Fermionic formulas and rigged configurations}
\label{sec:RC}
As before let $\la$ be a partition and $B=B^{r_k,s_k}\otimes\cdots
\otimes B^{r_1,s_1}$. Define the multiplicity array $L=(L_i^{(a)}\mid (a,i)\in \HH)$
where $L_i^{(a)}$ denotes the number of factors $B^{a,i}$ in $B$,
$\HH=I\times \Z_{>0}$ and $I=\{1,2,\ldots,n-1\}$.
The sequence of partitions $\nu=\{\nu^{(a)}\mid a\in I\}$
is an \textbf{$(L,\la)$-configuration} if
\begin{equation}\label{eq:constraint}
\sum_{(a,i)\in\HH} i m_i^{(a)} \alpha_a = \sum_{(a,i)\in\HH} i
L_i^{(a)} \La_a- \la,
\end{equation}
where $m_i^{(a)}$ is the number of parts of length $i$ in partition
$\nu^{(a)}$, $\La_a=\epsilon_1+\epsilon_2+\cdots+\epsilon_a$ are the
fundamental weights and $\alpha_a=\epsilon_a-\epsilon_{a+1}$ are the
simple roots of type $A_{n-1}$. Here $\epsilon_i$ is the $i$-th canonical
unit vector of $\Z^n$. The constraint \eqref{eq:constraint} is equivalent
to the condition
\begin{equation} \label{config def}
  |\nu^{(k)}| = \sum_{j>k} \la_j -
  	\sum_{a=1}^L s_a \max(r_a-k,0)
\end{equation}
on the size of $\nu^{(k)}$.

The \textbf{vacancy numbers} for the $(L,\la)$-configuration $\nu$
are defined as
\begin{equation*}
p_i^{(a)}=\sum_{j\ge 1} \min(i,j) L_j^{(a)}
 - \sum_{(b,j)\in \HH} (\alpha_a | \alpha_b) \min(i,j)m_j^{(b)},
\end{equation*}
where $(\cdot\mid\cdot)$ is the normalized invariant form on the weight lattice $P$
such that $(\alpha_a\mid \alpha_b)$ is the Cartan matrix.
The $(L,\la)$-configuration $\nu$ is \textbf{admissible} if $p^{(a)}_i\ge 0$ for 
all $(i,a)\in\HH$, and the set of admissible $(L,\la)$-configurations is denoted
by $\Confb(L,\la)$. 
It was proven in~\cite[Lemma 10]{KS:2002} that $p_i^{(a)}\ge 0$ for all existing 
parts $i$ implies that $p_i^{(a)}\ge 0$ for all $i$.

Set
\begin{equation*}
\cc(\nu)=\frac{1}{2} \sum_{a,b\in I} \sum_{j,k\ge 1} (\alpha_a \mid
\alpha_b) \min(j,k) m_j^{(a)} m_k^{(b)}.
\end{equation*}

With this notation we define the following fermionic formula.
It was first conjectured in~\cite{KS:2002,SW:1999} that it is an
explicit expression for the generalized Kostka polynomials,
stemming from the analogous expression of Kirillov and Reshetikhin~\cite{KR:1988} 
for the Kostka polynomial. This conjecture was proved in~\cite[Theorem 2.10]{KSS:2002}.

\begin{definition}[\textbf{Fermionic formula}]
For a multiplicity array $L$ and a partition $\la$ such that
$|\la|=\sum_{(a,i)\in\HH} aiL_i^{(a)}$ define
\begin{equation}\label{eq:fermi}
  \Mb(L,\la;q)=\sum_{\nu\in\Confb(L,\la)} q^{\cc(\nu)}
   \prod_{(a,i)\in\HH} \qbin{p_i^{(a)}+m_i^{(a)}}{m_i^{(a)}}.
\end{equation}
\end{definition}

Expression \eqref{eq:fermi} can be reformulated as the generating function
over rigged configurations. To this end we need
to define certain labelings of the rows of the
partitions in a configuration.
For this purpose one should view a partition as
a multiset of positive integers.
A rigged partition is by definition a finite multiset of
pairs $(i,x)$ where $i$ is a positive integer and
$x$ is a nonnegative integer.  The pairs $(i,x)$ are referred to
as strings; $i$ is referred to as the
length or size of the string and $x$ as the label or
quantum number of the string.  A rigged partition is
said to be a rigging of the partition $\rho$ if
the multiset, consisting of the sizes of the strings,
is the partition $\rho$.  So a rigging of $\rho$
is a labeling of the parts of $\rho$ by nonnegative integers,
where one identifies labelings that differ only by
permuting labels among equal sized parts of $\rho$.

A rigging $J$ of the $(L,\la)$-configuration $\nu$ is a sequence of
riggings of the partitions $\nu^{(a)}$ such that
every label $x$ of a part of $\nu^{(a)}$ of size $i$
satisfies the inequalities
\begin{equation} \label{eq:rigging def}
  0 \le x \le p^{(a)}_i.
\end{equation}
Alternatively, a rigging of a configuration $\nu$ may be viewed
as a double-sequence of partitions $J=(J^{(a,i)}\mid (a,i)\in\HH)$
where $J^{(a,i)}$ is a partition that has at most $m_i^{(a)}$ parts
each not exceeding $p_i^{(a)}$.
The pair $(\nu,J)$ is called a \textbf{rigged configuration}.
The set of riggings of admissible $(L,\la)$-configurations
is denoted by $\RCb(L,\la)$.
Let $(\nu,J)^{(a)}$ be the $a$-th rigged partition
of $(\nu,J)$.  A string $(i,x)\in (\nu,J)^{(a)}$
is said to be \textbf{singular} if $x=p^{(a)}_i$, that is,
its label takes on the maximum value.

\begin{example}\label{ex:rc}
Let $L$ be the multiplicity array of $B=(B^{1,1})^{\otimes 2} \otimes B^{1,4}
\otimes B^{2,1} \otimes B^{2,3}$ and $\la=(6,4,3,1)$. Then
\begin{equation*}
(\nu,J)= \yngrc(3,{1,1},1,{1,1}) \quad \yngrc(3,{1,1},1,{0,1}) \quad 
\yngrc(1,{0,0})\in \RCb(L,\la),
\end{equation*}
where the first number next to each part is the rigging and the second
one is the vacancy number for the corresponding part.
\end{example}

The set of rigged configurations is endowed with a natural
statistic $\cc$ defined by
\begin{equation} \label{eq:RC charge}
  \cc(\nu,J)=\cc(\nu)+\sum_{(a,i)\in\HH} |J^{(a,i)}|
\end{equation}
for $(\nu,J)\in\RCb(L,\la)$, where $|J^{(a,i)}|$ is the size of partition
$J^{(a,i)}$.
Since the $q$-binomial $\qbins{p+m}{m}$ is the generating
function of partitions with at most $m$ parts each not
exceeding $p$, \eqref{eq:fermi} can be rewritten as
\begin{equation}\label{eq:rc}
  \Mb(L,\la;q)=\sum_{(\nu,J)\in\RCb(L,\la)} q^{\cc(\nu,J)}.
\end{equation}

The $\Xb=\Mb$ conjecture asserts that $\Mb(L,\la;q)=\Xb(B,\la;q)$ where
$L$ is the multiplicity array of $B$.
For type $A$ this was proven in~\cite{KSS:2002} by showing that there is
a bijection $\Phib:\Pathb(B,\la)\to\RCb(L,\la)$ which preserves the statisitics.
\begin{theorem}\cite[Theorem 2.12]{KSS:2002}\label{thm:X=M}
For $\la$ a partition, $B^{r_k,s_k}\otimes B^{r_{k-1},s_{k-1}}\otimes \cdots 
\otimes B^{r_1,s_1}$ and $L$ the corresponding multiplicity array such that
$|\la|=\sum_j r_js_j$ we have $\Mb(L,\la;q)=\Xb(B,\la;q)$.
\end{theorem}

\subsection{Operations on crystals}
\label{sec:operations crystals}

To define the bijection $\Phib$ we first need to define certain maps on paths and
rigged configurations. These maps correspond to the following operations on crystals:
\begin{enumerate}
\item If $B=B^{1,1}\otimes B'$, let $\lh(B)=B'$. This operation is called 
\textbf{left-hat}.
\item If $B=B^{r,s}\otimes B'$ with $s\ge 2$, let $\ls(B)=B^{r,1}\otimes
 B^{r,s-1}\otimes B'$. This operation is called \textbf{left-split}.
\item If $B=B^{r,1}\otimes B'$ with $r\ge 2$, let $\lb(B)=B^{1,1}\otimes B^{r-1,1}
\otimes B'$. This operation is called \textbf{box-split}.
\end{enumerate}
In analogy we define $\lh(L)$ (resp. $\ls(L)$, $\lb(L)$) to be the multiplicity array of 
$\lh(B)$ (resp. $\ls(B)$, $\lb(B)$), if $L$ is the multiplicity array of $B$.
The corresponding maps on crystal elements are given by:
\begin{enumerate}
\item Let $b=c\otimes b'\in B^{1,1}\otimes B'$. Then $\lh(b)=b'$.
\item Let $b=c\otimes b'\in B^{r,s}\otimes B'$, where $c=c_1c_2\cdots c_s$ and $c_i$
denotes the $i$-th column of $c$. Then $\ls(b)=c_1\otimes c_2\cdots c_s\otimes b'$.
\item Let $b=\begin{array}{|c|} \hline b_1\\ \hline b_2\\ \hline \vdots\\ \hline b_r\\ 
\hline \end{array}\otimes b'\in B^{r,1}\otimes B'$, where $b_1<\cdots<b_r$.
Then $\lb(b)=\begin{array}{|c|} \hline b_r\\ \hline \end{array} \otimes 
\begin{array}{|c|} \hline b_1\\ \hline \vdots \\ \hline b_{r-1}\\ \hline \end{array} 
\otimes b'$.
\end{enumerate}

In the next subsection we define the corresponding maps on rigged configurations,
and give the bijection in subsection~\ref{sec:bij}.

\subsection{Operations on rigged configurations}
\label{sec:operations rc}
Suppose $L_1^{(1)}>0$. The main algorithm on rigged configurations as defined
in~\cite{KR:1988,KSS:2002} for admissible rigged configurations is called $\delta$.
For a partition $\la=(\la_1,\ldots,\la_n)$,
let $\lm$ be the set of all nonnegative tuples $\mu=(\mu_1,\ldots,\mu_n)$ such that
$\la-\mu=\epsilon_r$ for some $1\le r\le n$.
Define $\delta:\RCb(L,\la)\to \bigcup_{\mu\in\lm} \RCb(\lh(L),\mu)$
by the following algorithm. Let $(\nu,J)\in\RCb(L,\la)$. Set $\ell^{(0)}=1$ and repeat the
following process for $a=1,2,\ldots,n-1$ or until stopped. Find the smallest index 
$i\ge \ell^{(a-1)}$ such that $J^{(a,i)}$ is singular. If no such $i$ exists, set 
$\rk(\nu,J)=a$ and stop. Otherwise set $\ell^{(a)}=i$ and continue with $a+1$.
Set all undefined $\ell^{(a)}$ to $\infty$.

The new rigged configuration $(\tilde{\nu},\tilde{J})=\delta(\nu,J)$ is obtained by
removing a box from the selected strings and making the new strings singular
again. Explicitly
\begin{equation*}
 m_i^{(a)}(\tilde{\nu})=m_i^{(a)}(\nu)+\begin{cases}
 1 & \text{if $i=\ell^{(a)}-1$}\\
 -1 & \text{if $i=\ell^{(a)}$}\\
 0 & \text{otherwise.} \end{cases}
\end{equation*}
The partition $\tilde{J}^{(a,i)}$ is obtained from $J^{(a,i)}$ by removing
a part of size $p_i^{(a)}(\nu)$ for $i=\ell^{(a)}$,
adding a part of size $p_i^{(a)}(\tilde{\nu})$ for $i=\ell^{(a)}-1$, 
and leaving it unchanged otherwise. Then $\delta(\nu,J)\in \RCb(\lh(L),\mu)$
where $\mu=\la-\epsilon_{\rk(\nu,J)}$.

\begin{example}\label{ex:delta}
Let $(\nu,J)$ be the rigged configuration of Example~\ref{ex:rc}.
Hence $\ell^{(1)}=1$, $\ell^{(2)}=3$ and $\ell^{(3)}=\infty$, so that
$\rk(\nu,J)=3$ and 
\begin{equation*}
\delta(\nu,J)= \yngrc(3,{1,1}) \quad \yngrc(2,{0,0},1,{0,0}) \quad \yngrc(1,{0,0}).
\end{equation*}
Also $\cc(\nu,J)=8$.
\end{example}

Let $s\ge2$. Suppose $B=B^{r,s}\otimes B'$ and $L$ the corresponding 
multiplicity array. Note that $\Confb(L,\la)\subset \Confb(\ls(L),\la)$. Under this
inclusion map, the vacancy number $p_i^{(a)}$ for $\nu$ increases by
$\delta_{a,r} \chi(i<s)$. Hence there is a well-defined injective map
$\rcls:\RCb(L,\la)\rightarrow \RCb(\ls(L),\la)$ given by $\rcls(\nu,J)=(\nu,J)$.

Suppose $r\ge2$ and $B=B^{r,1}\otimes B'$ with multiplicity array $L$.
Then there is an injection $\rclb:\RCb(L,\la)\to \RCb(\lb(L),\la)$ defined by adding 
singular strings of length $1$ to $(\nu,J)^{(a)}$ for $1\le a < r$. Moreover the
vacancy numbers stay the same.

\subsection{Bijection}\label{sec:bij}
The map $\Phib:\Pathb(B,\la)\to\RCb(L,\la)$ is defined by various commutative diagrams.
Note that it is possible to go from $B=B^{r_k,s_k}\otimes B^{r_{k-1},s_{k-1}}\otimes \cdots
\otimes B^{r_1,s_1}$ to the empty crystal via successive application of $\lh$, $\ls$ and
$\lb$.

\begin{definition} \label{def:bij}
Define that map $\Phib:\Pathb(B,\la)\rightarrow \RCb(L,\la)$ such that 
the empty path maps to the empty rigged configuration, and:
\begin{enumerate}
\item Suppose $B=B^{1,1} \otimes B'$. Then the diagram
\begin{equation*}
\begin{CD}
\Pathb(B,\la) @>{\phi}>> \RCb(L,\la) \\
@V{\lh}VV @VV{\delta}V \\
\displaystyle{\bigcup_{\mu\in\lm} \Pathb(\lh(B),\mu)} @>>{\phi}> \displaystyle{\bigcup_{\mu\in\lm}
\RCb(\lh(L),\mu)}
\end{CD}
\end{equation*}
commutes.
\item Suppose $B=B^{r,s} \otimes B'$ with $s\ge 2$. Then the following diagram commutes:
\begin{equation*}
\begin{CD}
\Pathb(B,\la) @>{\phi}>> \RCb(L,\la) \\
@V{\ls}VV @VV{\rcls}V \\
\Pathb(\ls(B),\la) @>>{\phi}> \RCb(\ls(L),\la)
\end{CD}
\end{equation*}
\item Suppose $B=B^{r,1} \otimes B'$ with $r\ge2$. Then the following diagram commutes:
\begin{equation*}
\begin{CD}
\Pathb(B,\la) @>{\phi}>> \RCb(L,\la) \\
@V{\lb}VV @VV{\rclb}V \\
\Pathb(\lb(B),\la) @>>{\phi}> \RCb(\lb(L),\la)
\end{CD}
\end{equation*}
\end{enumerate}
\end{definition}

\begin{theorem}\cite{KSS:2002}\label{thm:bij}
The map $\Phib:\Pathb(B,\la)\to\RCb(L,\la)$ is a bijection and preserves the statistics,
that is, $\Dt(b)=\cc(\Phib(b))$ for all $b\in\Pathb(B,\la)$.
\end{theorem}
Note that Theorem~\ref{thm:bij} immediately implies Theorem~\ref{thm:X=M}.

\begin{example}\label{ex:path}
The path which corresponds to $(\nu,J)$ of Example~\ref{ex:rc} under $\Phib$ is
\begin{equation*}
b=\young(3) \otimes \young(2) \otimes \young(1134) \otimes \young(1,3) \otimes
\young(111,222) \in\Pathb(B,\la).
\end{equation*}
We have $\Dt(b)=\cc(\nu,J)=8$. The steps of Definition~\ref{def:bij} are summarized
in Table~\ref{tab:ex path}.
\begin{table}
\begin{equation*}
\begin{array}{|c|c|r|}
\hline
\mathrm{Step} & (\nu,J) & b \phantom{bbbbbbbbbbbbbbbbbbb}\\[2mm] \hline
& \yngrc(3,{1,1},1,{1,1}) \quad \yngrc(3,{1,1},1,{0,1}) \quad 
\yngrc(1,{0,0})
& \young(3) \otimes \young(2) \otimes \young(1134) \otimes \young(1,3) \otimes
\young(111,222)\\[3mm]
(1) & \yngrc(3,{1,1}) \quad \yngrc(2,{0,0},1,{0,0}) \quad \yngrc(1,{0,0})
& \young(2) \otimes \young(1134) \otimes \young(1,3) \otimes
\young(111,222)\\[3mm]
(1) & \yngrc(2,{1,1}) \quad \yngrc(2,{0,0},1,{0,0}) \quad \yngrc(1,{0,0})
& \young(1134) \otimes \young(1,3) \otimes \young(111,222)\\[3mm]
(2) & \yngrc(2,{1,2}) \quad \yngrc(2,{0,0},1,{0,0}) \quad \yngrc(1,{0,0})
& \young(1) \otimes \young(134) \otimes \young(1,3) \otimes \young(111,222)\\[3mm]
(1) & \yngrc(2,{1,1}) \quad \yngrc(2,{0,0},1,{0,0}) \quad \yngrc(1,{0,0})
& \young(134) \otimes \young(1,3) \otimes \young(111,222)\\[3mm]
(2) & \yngrc(2,{1,2}) \quad \yngrc(2,{0,0},1,{0,0}) \quad \yngrc(1,{0,0})
& \young(1) \otimes \young(34) \otimes \young(1,3) \otimes \young(111,222)\\[3mm]
(1) & \yngrc(2,{1,1}) \quad \yngrc(2,{0,0},1,{0,0}) \quad \yngrc(1,{0,0})
& \young(34) \otimes \young(1,3) \otimes \young(111,222)\\[3mm]
(2) & \yngrc(2,{1,1}) \quad \yngrc(2,{0,0},1,{0,0}) \quad \yngrc(1,{0,0})
& \young(3) \otimes \young(4) \otimes \young(1,3) \otimes \young(111,222)\\[3mm]
(1) & \yngrc(1,{1,1}) \quad \yngrc(1,{0,0},1,{0,0}) \quad \yngrc(1,{0,0})
& \young(4) \otimes \young(1,3) \otimes \young(111,222)\\[3mm]
(1) & \emptyset \quad \yngrc(1,{0,0}) \quad \emptyset
& \young(1,3) \otimes \young(111,222)\\[3mm]
(3) & \yngrc(1,{1,1}) \quad \yngrc(1,{0,0}) \quad \emptyset
& \young(3) \otimes \young(1) \otimes \young(111,222)\\[3mm]
(1) & \emptyset \quad \emptyset \quad \emptyset
& \young(1) \otimes \young(111,222)\\[3mm]
(1) & \emptyset \quad \emptyset \quad \emptyset
& \young(111,222)\\[3mm]
\hline
\end{array}
\end{equation*}
\caption{\label{tab:ex path} Explicit steps for Example~\ref{ex:path}}
\end{table}
\end{example}

\subsection{Properties}
\label{sec:prop}
As we have already seen in Section~\ref{sec:bij}, the bijection $\Phib$
preserves the statistics. In addition to this it satisfies a couple of
other amazing properties, one of them being the evacuation theorem.
The Dynkin diagram of type $A_{n-1}$ has the symmetry $\tau$ which
interchanges $i$ and $n-i$. There is a corresponding map $*$ on crystals
which satisfies
\begin{equation}\label{eq:*}
\begin{split}
\wt(b^*)&= w_0 \wt(b)\\
e_i(b^*)&= f_{\tau(i)}(b)^* \\
f_i(b^*)&= e_{\tau(i)}(b)^*
\end{split}
\end{equation}
for all $i\in I$ where $w_0$ is the longest permutation of the symmetric group
$S_{n-1}$. Explicitly an element $i\in B^{1,1}$ is mapped to $n+1-i$.
For $b\in B^{r,s}$, $b^*$ is the tableau obtained by replacing every entry $c$ 
of $b$ by $c^*$ and then rotating by 180 degrees. The resulting tableau is sometimes
called the \textbf{antitableau} of $b$.
For $b=b_k\otimes b_{k-1}\otimes \cdots\otimes b_1\in B^{r_k,s_k}\otimes 
B^{r_{k-1},s_{k-1}}\otimes \cdots \otimes B^{r_1,s_1}$ define
$b^*=b_1^*\otimes b_2^*\otimes \cdots \otimes b_k^*$.
\begin{example}
For type $A_4^{(1)}$
\begin{equation*}
{\young(11,23)}^{\; *}=\young(34,55)\,.
\end{equation*}
\end{example}
By \eqref{eq:*} the map $*$ maps classical components to classical components.
By weight considerations, these components have to be of the same classical
highest weight. Let $\ev(b)$ be the highest weight vector in the same
classical component as $b^*$.

\begin{example}
Let $b$ be the path of Example \ref{ex:path}. Then
\begin{equation*}
b^*= \young(333,444) \otimes \young(2,4) \otimes \young(1244) \otimes 
\young(3) \otimes \young(2)
\end{equation*}
and
\begin{equation*}
\ev(b)=\young(112,234) \otimes \young(2,3) \otimes \young(1113) \otimes \young(2)
\otimes \young(1).
\end{equation*}
\end{example}

On rigged configurations define $\theta$ to be the complementation of
quantum numbers. More precisely, if $(i,x)$ is a string in $(\nu,J)^{(k)}$,
replace this string by $(i,p_i^{(k)}-x)$. The 
\textbf{Evacuation Theorem}~\cite[Theorem 5.6]{KSS:2002} asserts that $\ev$ and 
$\theta$ correspond under the bijection $\Phib$.
\begin{example}
For $(\nu,J)$ of Example~\ref{ex:rc} we have
\begin{equation*}
\theta(\nu,J)= \yngrc(3,{0,1},1,{0,1}) \quad \yngrc(3,{0,1},1,{1,1}) \quad 
\yngrc(1,{0,0})
\end{equation*}
and it is easy to check that $\theta(\Phib(b))=\Phib(\ev(b))$ with $b$ as in
Example~\ref{ex:path}.
\end{example}

The combinatorial $R$ matrix on crystals is the identity on rigged configurations
under the bijection $\Phib$. See for example \cite[Lemma 8.5]{KSS:2002}
or~\cite[Theorem 8.6]{SS:2004}. This shows in particular that the polynomial
$\Xb(B,\la;q)$ does not depend on the order of the tensor factors in $B$.

\begin{example}
Take $b$ from Example~\ref{ex:path}. Then $R_1$ is the combinatorial $R$-matrix
applied to the first two tensor factors and
\begin{equation*}
R_1(b)=\young(3) \otimes \young(2) \otimes \young(1134) \otimes 
\young(111,223) \otimes \young(1,2).
\end{equation*}
It can be checked that $\Phib(b)=\Phib(R_1(b))=(\nu,J)$ is the rigged configuration
of Example~\ref{ex:rc}.
\end{example}

The bijection $\Phib$ is also well-behaved with respect to \textbf{transpose duality}. 
Define
\begin{equation*}
\trP:B=B^{r_k,s_k}\otimes B^{r_{k-1},s_{k-1}} \otimes \cdots \otimes B^{r_1,s_1} \to 
B^{s_k,r_k}\otimes B^{s_{k-1},r_{k-1}} \otimes \cdots \otimes B^{s_1,r_1}=:B^t
\end{equation*}
as follows. For $b=b_k\otimes \cdots \otimes b_1\in
B^{r_k,s_k}\otimes \cdots \otimes B^{r_1,s_1}$ rotate each rectangular tableau
$b_i$ by $90^o$ clockwise to obtain $\tilde{b}_i$. Suppose the letter $a$ occurs
in cell $c$ of $\tilde{b}_i$. Then replace letter $a$ in cell $c$ by $\tilde{a}$
where $\tilde{a}$ is chosen such that the letter $a$ in cell $c$ is the 
$\tilde{a}$-th letter $a$ in $\row(b)$ reading from right to left.
Since heighest-weight crystal elements are mapped to heighest-weight elements
this induces a map 
\begin{equation*}
\trP:\Pathb(B,\la) \to \Pathb(B^t,\la^t).
\end{equation*}
It should be noted that we are assuming here that $n$ is big enough so that both
$B^{r_i,s_i}$ and $B^{s_i,r_i}$ are $A_{n-1}^{(1)}$ crystals.

The analogous map on rigged configurations is
\begin{equation*}
\trRC:\RCb(L,\la)\to\RCb(L^t,\la^t),
\end{equation*}
where $L^t$ is the multiplicity array of $B^t$.
Let $(\nu,J)\in\RCb(L,\la)$ and let $\nu$ have
the associated matrix $m$ with entries $m_{ai}$ as in \cite[(9.2)]{KS:2002}
\begin{equation*}
  m_{ai} = \sum_{i\le j} m_j^{(a-1)} - m_j^{(a)}.
\end{equation*}
Note that $\sum_{i\le j} m_j^{(a)}$ is the size of the
$i$-th column of the partition $\nu^{(a)}$. Here $m_j^{(0)}$ is defined to be zero.
The configuration $\nu^t$ in $(\nu^t,J^t)=\trRC(\nu,J)$ is defined
by its associated matrix $\mt$ given by
\begin{equation*}
  \mt_{ai} = - m_{ia} +\chi((i,a)\in \la) -
	\sum_{j=1}^k \chi((i,a)\in (s_j^{r_j})).
\end{equation*}
Here $(i,a) \in \la$ means that
the cell $(i,a)$ is in the Ferrers diagram of the partition $\la$
with $i$ specifying the row and $a$ the column.

Recall that the riggings $J$ can be viewed as a double sequence of 
partitions $J=(J^{(a,i)})$ where $J^{(a,i)}$ is a partition inside the
rectangle of height $m_i^{(a)}$ and width $p_i^{(a)}$.
The partition $J^{t(i,a)}$ corresponding to $(\nu^t,J^t)=\trRC(\nu,J)$ is 
defined as the transpose of the complementary partition to $J^{(a,i)}$ in 
the rectangle of height $m_i^{(a)}$ and width $p_i^{(a)}$.

The \textbf{Transpose Theorem} \cite[Theorem 7.1]{KSS:2002} asserts that
$\Phib(\trP(b))=\trRC(\Phib(b))$ for all $b\in\Pathb(B,\la)$. This implies 
in particular the transpose symmetry \cite[Theorem 7.1]{SW:1999} 
\cite[Conjecture 3]{KS:2002}
\begin{equation*}
\Xb(B,\la;q)=\widetilde{\Xb}(B^t,\la^t;q)
\end{equation*}
and similarly for $\Mb$, where $\widetilde{\Xb}(B,\la;q)=q^{n(B)}\Xb(B,\la;q^{-1})$
and
\begin{equation*}
n(B)=\sum_{1\le i<j\le k} \min(s_i,s_j)\min(r_i,r_j).
\end{equation*}

\begin{example}
As usual let $b$ be the path of Example~\ref{ex:path}. Then
\begin{equation*}
\trP(b)= \young(3) \otimes \young(4) \otimes \young(1,2,5,6) \otimes \young(14)
\otimes \young(11,22,33).
\end{equation*}
Similarly, let $(\nu,J)$ be the rigged configuration of Example~\ref{ex:rc}. 
Then the matrix $m$ and $m^t$ are
\begin{equation*}
\begin{split}
m&=\begin{pmatrix} 2 & -1 & -1 & \cdots \\ 0 & 0 & 0 & \cdots \\ 1 & 1 & 1 & 
\cdots\\ 1 & 0 & 0 & \cdots\\ \vdots & \vdots &\vdots& \end{pmatrix}\\
m^t&=\begin{pmatrix} -2 & -1 & \cdots\\ 0&0&\cdots\\ 0&0&\cdots\\
0&1&\cdots\\ 1&0&\cdots \\1&0&\cdots \\ \vdots&\vdots& \end{pmatrix}
\end{split}
\end{equation*}
so that
\begin{equation*}
\trRC(\nu,J)= \yngrc(2,{1,1},1,{0,1}) \quad \yngrc(2,{0,0},1,{0,0}) \quad
\yngrc(2,{0,1},1,{0,1}) \quad \yngrc(1,{0,0},1,{0,0}) \quad \yngrc(1,{0,0}).
\end{equation*}
It can be checked explicitly in this example that $\Phib(\trP(b))=\trRC(\Phib(b))$.
\end{example}

Finally let us mention the \textbf{contragredient duality} which is of great
importance for the notion of virtual crystals~\cite{OSS:2003a,OSS:2003b}.
On crystals define the map
\begin{equation*}
\vee : B^{r,s}\to B^{n-r,s}
\end{equation*}
where each column $c=c_1\ldots c_r$ of $b\in B^{r,s}$ is replaced by
column $(n+1-d_{n-r}) \ldots (n+1-d_1)$ where $\{d_1<d_2<\cdots<d_{n-r}\}$
is the complement of $\{c_1<c_2<\ldots<c_r\}$ in $\{1,2,\ldots,n\}$. 
Note that $e_i(b)^\vee=f_{n-i}(b^\vee)$.
\begin{example}
The contragredient dual of $b=\young(12,23)$ for $n=4$ is 
$b^\vee=\young(11,24)$.
\end{example}
The map $\vee$ can be extended to a map on paths
\begin{equation*}
\vee:\Pathb(B,\la)\to \Pathb(B^\vee,\la^\vee),
\end{equation*}
where $B^\vee=B^{n-r_k,s_k}\otimes \cdots\otimes B^{n-r_1,s_1}$,
$\la^\vee=(N-\la_n,N-\la_{n-1},\ldots,N-\la_1)$ and $N=s_1+\cdots+s_k$,
by mapping $b=b_k\otimes \cdots \otimes b_1$ to 
$b^\vee=b_k^\vee \otimes \cdots \otimes b_1^\vee$.

For given $n$, define 
\begin{equation*}
\rev:\RCb(L,\la)\to \RCb(L^\vee,\la^{\vee})
\end{equation*}
such that for $(\nu^\vee,J^\vee)=\rev(\nu,J)$ we have
$(\nu^\vee,J^\vee)^{(a)}=(\nu,J)^{(n-a)}$. 
Then we have~\cite[Theorem 5.7]{OSS:2003a} that
$\Phib(b^\vee)=\rev(\Phib(b))$ for all $b\in\Pathb(B,\la)$.
This implies the contragredient symmetry
\begin{equation*}
\Xb(B^\vee,\la^\vee;q)=\Xb(B,\la;q)
\end{equation*}
and similarly for $\Mb$.

\begin{example}
Employing one last time $b$ of Example~\ref{ex:path} we obtain
\begin{equation*}
b^\vee= \young(1,3,4) \otimes \young(1,2,4) \otimes \young(1112,2233,3344)
\otimes \young(1,3) \otimes \young(111,222)
\end{equation*}
and for $(\nu,J)$ of Example~\ref{ex:rc}
\begin{equation*}
\rev(\nu,J)=\yngrc(1,{0,0}) \quad \yngrc(3,{1,1},1,{0,1}) \quad
\yngrc(3,{1,1},1,{1,1})
\end{equation*}
which is also $\Phib(b^\vee)$.
\end{example}

The bijection $\Phib$ has further properties. For example it is well-behaved
under certain embeddings. We refer the interested reader to the 
literature~\cite{KSS:2002,SW:1999,KS:2002,SS:2004}.

\subsection{Open Problems}
\begin{itemize}
\item For nonexceptional types, the bijection $\Phib$ was given in~\cite{OSS:2003,SS:2004}
for the cases $B=B^{1,s_k}\otimes \cdots \otimes B^{1,s_1}$ and for type $D_n^{(1)}$
in the case $B=B^{r_k,1}\otimes \cdots \otimes B^{r_1,1}$~\cite{S:2005}. For all other cases,
it is still an outstanding problem to prove that $\Phib$ exists. In particular, the analogues
of the splitting maps need to be found.
\item It would be very nice to have a more conceptual definition of the bijection $\Phib$
rather than the recursive definition in terms of the splitting and hatting maps.
A possible avenue would be to give a definition of $\Phib$ in terms of the affine crystal
structure on rigged configurations. In section~\ref{sec:unrestricted} we provide such
a crystal structure for $B^{r,s}$ of type $A_{n-1}^{(1)}$. To obtain $\Phib$,
one would need the affine crystal structure on tensor products 
$B=B^{r_k,s_k}\otimes \cdots \otimes B^{r_1,s_1}$. Compare with section~\ref{sec:open X=M}.
\end{itemize}

\section{$X=M$}
\label{sec:unrestricted}

In this section we deal with the unrestricted version of the $X=M$
conjecture for type $A_{n-1}^{(1)}$. In particular it is our aim to find
a fermionic formula for the unrestricted configuration sum $X(B,\lambda;q)$
of Definition~\ref{def:X}. This has recently been achieved in~\cite{S:2005a}
by extending the set of rigged configurations to the set of unrestricted
rigged configurations by imposing a crystal structure in this set.
A direct bijection between unrestricted paths and unrestricted rigged configurations
along the lines of Definition~\ref{def:bij} was given in~\cite{DS:2005}.
Here we mostly follow~\cite{S:2005a} and derive the fermionic formula
$M(B,\lambda;q)$ from the crystal structure on rigged configurations.

\subsection{Crystal structure on rigged configurations}\label{sec:rc crystal}
The set of unrestricted rigged configurations $\RC(L)$ can be introduced by defining 
a crystal structure generated from highest weight vectors 
given by elements in $\RCb(L)=\bigcup_{\la} \RCb(L,\la)$ by the Kashiwara 
operators $e_a,f_a$.

\begin{definition} \label{def:crystal} Let $L$ be a multiplicity array.
Define the set of \textbf{unrestricted rigged configurations} $\RC(L)$
as the set generated from the elements in $\RCb(L)$ by the application of  
the operators $f_a,e_a$ for $a\in I$ defined as follows:
\begin{enumerate}
\item
Define $e_a(\nu,J)$ by removing a box from a string of length $k$ in
$(\nu,J)^{(a)}$ leaving all colabels fixed and increasing the new
label by one. Here $k$ is the length of the string with the smallest
negative rigging of smallest length. If no such string exists,
$e_a(\nu,J)$ is undefined. 
\item
Define $f_a(\nu,J)$ by adding a box to a string of length $k$ in
$(\nu,J)^{(a)}$ leaving all colabels fixed and decreasing the new
label by one. Here $k$ is the length of the string with the smallest
nonpositive rigging of largest length. If no such string exists,
add a new string of length one and label -1.
If the result is not a valid unrestricted rigged configuration
$f_a(\nu,J)$ is undefined.
\end{enumerate}
\end{definition}

Let $(\nu,J)\in \RC(L)$. If $f_a$ adds a box to a string of length $k$ in
$(\nu,J)^{(a)}$, then the vacancy numbers change according to
\begin{equation}\label{eq:change in p}
p_i^{(b)} \mapsto p_i^{(b)}-(\alpha_a | \alpha_b)\chi(i>k),
\end{equation}
where $\chi(S)=1$ if the statement $S$ is true and $\chi(S)=0$ if $S$
is false. Similarly, if $e_a$ adds a box of length $k$ to $(\nu,J)^{(a)}$,
then the vacancy numbers change as
\begin{equation*}
p_i^{(b)} \mapsto p_i^{(b)}+(\alpha_a | \alpha_b)\chi(i\ge k).
\end{equation*}

We may define a weight function $\wt:\RC(L)\to P$ as
\begin{equation}\label{eq:rc weight}
\wt(\nu,J)=\sum_{(a,i)\in \HH} i(L_i^{(a)}\La_a - m_i^{(a)}\alpha_a)
\end{equation}
for $(\nu,J)\in\RC(L)$. It is clear from the definition that 
$\wt(f_a (\nu,J))=\wt(\nu,J)-\alpha_a$.
Define
\begin{equation*}
\RC(L,\la)=\{(\nu,J)\in \RC(L) \mid \wt(\nu,J)=\la\}.
\end{equation*}

\begin{example} Let $\mathfrak{g}$ be of type $A_2^{(1)}$. Let $\la=(3,2,3)$,
$L_1^{(1)}=L_3^{(1)}=L_2^{(2)}=1$ and all other $L_i^{(a)}=0$. Then
\begin{equation*}
(\nu,J)= \quad \yngrc(2,-1,1,-1) \quad \yngrc(3,-2)
\end{equation*}
is in $\RC(L,\la)$, where the parts of the rigging $J^{(a,i)}$ are written
next to the parts of length $i$ in partition $\nu^{(a)}$. We have
\begin{equation*}
f_1(\nu,J)= \yngrc(3,-2,1,-1) \quad \yngrc(3,-1) \quad \text{and} \quad
e_1(\nu,J)= \yngrc(2,1) \quad \yngrc(3,-3).
\end{equation*}
\end{example}

\begin{example} Let $\mathfrak{g}$ be of type $A_2^{(1)}$. Let $\la=(4,5,6)$,
$L_1^{(1)}=15$ and all other $L_i^{(a)}=0$. Then
\begin{equation*}
(\nu,J)= \quad \yngrc(3,-2,2,0,2,-1,1,2,1,2,1,0,1,0) \quad
\yngrc(3,-1,2,-1,1,1)
\end{equation*}
is in $\RC(L,\la)$. We have
\begin{equation*}
e_1(\nu,J)= \yngrc(2,0,2,-1,2,-1,1,2,1,2,1,0,1,0) \quad
\yngrc(3,-2,2,-1,1,1) \quad \text{and} \quad
e_2(\nu,J)= \yngrc(3,-3,2,-1,2,-2,1,2,1,2,1,0,1,0) \quad
\yngrc(3,1,1,1,1,0).
\end{equation*}
\end{example}

The following Theorem was proven in~\cite{S:2005a} for all simply-laced algebras.
\begin{theorem}\cite[Theorem 3.7]{S:2005a} \label{thm:rc crystal}
The graph generated from $(\nub,\Jb)\in \RCb(L,\la)$ and the crystal operators
$e_a,f_a$ of Definition~\ref{def:crystal} is isomorphic to the crystal graph $B(\la)$
of highest weight $\la$.
\end{theorem}

\begin{example} \label{ex:B21}
Consider the crystal $B(\yng(2,1))$ of type $A_2$
in $B=(B^{1,1})^{\otimes 3}$. Here is the crystal graph in the usual
labeling and the rigged configuration labeling:

\begin{picture}(100,200)(0,0)
\put(45,180){121}
\put(10,140){221}
\put(10,100){231}
\put(10,60){331}
\put(45,20){332}
\put(80,140){131}
\put(80,100){132}
\put(80,60){232}
\LongArrow(50,178)(25,152)
\LongArrow(60,178)(85,152)
\LongArrow(20,138)(20,112)
\LongArrow(20,98)(20,72)
\LongArrow(90,138)(90,112)
\LongArrow(90,98)(90,72)
\LongArrow(25,58)(50,32)
\LongArrow(85,58)(60,32)
\PText(30,165)(0)[b]{1}
\PText(15,120)(0)[b]{2}
\PText(15,80)(0)[b]{2}
\PText(30,40)(0)[b]{1}
\PText(95,120)(0)[b]{1}
\PText(95,80)(0)[b]{1}
\PText(80,40)(0)[b]{2}
\PText(80,165)(0)[b]{2}
\end{picture}
\hspace{1.5cm}
\begin{picture}(100,200)(0,0)
\put(40,185){{\tiny $\yngrc(1,0) \; \emptyset$}}
\put(0,145){{\tiny $\yngrc(2,-1) \; \emptyset$}}
\put(0,105){{\tiny $\yngrc(2,0) \; \yngrc(1,-1)$}}
\put(0,65){{\tiny $\yngrc(2,1) \; \yngrc(2,-2)$}}
\put(30,17){{\tiny $\yngrc(2,-1,1,-1) \; \yngrc(2,-1)$}}
\put(80,145){{\tiny $\yngrc(1,1) \; \yngrc(1,-1)$}}
\put(80,103){{\tiny $\yngrc(1,-1,1,-1) \; \yngrc(1,0)$}}
\put(80,63){{\tiny $\yngrc(2,-2,1,-1) \; \yngrc(1,0)$}}
\LongArrow(50,176)(25,152)
\LongArrow(60,176)(85,152)
\LongArrow(20,138)(20,112)
\LongArrow(20,97)(20,72)
\LongArrow(90,138)(90,117)
\LongArrow(90,94)(90,77)
\LongArrow(25,56)(50,32)
\LongArrow(83,53)(60,32)
\PText(30,165)(0)[b]{1}
\PText(15,120)(0)[b]{2}
\PText(15,80)(0)[b]{2}
\PText(30,40)(0)[b]{1}
\PText(95,120)(0)[b]{1}
\PText(95,80)(0)[b]{1}
\PText(80,40)(0)[b]{2}
\PText(80,165)(0)[b]{2}
\end{picture}
\end{example}

\begin{theorem} \cite[Theorem 3.9]{S:2005a}\label{thm:charge}
The cocharge $\cc$ as defined in~\eqref{eq:RC charge} is constant on connected
crystal components.
\end{theorem}

\begin{example}
The cocharge of the connected component in Example~\ref{ex:B21} is 1.
\end{example}

Combining the various results yields a generalization of Theorem~\ref{thm:bij}.
\begin{theorem} \cite[Theoren 3.10]{S:2005a}\label{thm:bij new}
Let $\la$ be a composition, $B$ be as in Theorem~\ref{thm:bij} and $L$ the corresponding
multiplicity array. Then there is a bijection $\Phi:\Path(B,\la)\to\RC(L,\la)$ 
which preserves the statistics, that is, $\Dt(b)=\cc(\Phi(b))$ for all $b\in\Path(B,\la)$.
\end{theorem}
\begin{proof}
By Theorem~\ref{thm:bij} there is such a bijection for the maximal elements
$b\in \Pathb(B)$. By Theorems~\ref{thm:rc crystal} and~\ref{thm:charge} this 
extends to all of $\Path(B,\la)$.
\end{proof}

Extending the definition of~\eqref{eq:rc} to
\begin{equation}\label{eq:XM def}
M(L,\la;q)=\sum_{(\nu,J)\in \RC(L,\la)} q^{\cc(\nu,J)},
\end{equation}
we obtain the corollary:
\begin{corollary} \cite[Corollary 3.10]{S:2005a}\label{cor:X=M}
With all hypotheses of Theorem~\ref{thm:bij new}, we have $X(B,\la;q)=M(L,\la;q)$.
\end{corollary}

\begin{example}
Let $n=4$, $B=B^{2,2}\otimes B^{2,1}$ and $\la=(2,2,1,1)$. Then the multiplicity
array is $L_1^{(2)}=1,L_2^{(2)}=1$ and $L_i^{(a)}=0$ for all other $(a,i)$. There are 
7 possible unrestricted paths in $\Path(B,\la)$. For each path $b\in\Path(B,\la)$ the 
corresponding rigged configuration $(\nu,J)=\Phi(b)$ together with the tail energy 
and cocharge is summarized below.
\begin{equation*}
\begin{array}{llllll}
b\;=\; \young(11,22) \otimes \young(3,4)
&
\quad (\nu,J) \;=
& \yngrc(1,0) & \yngrc(1,-1,1,-1) & \yngrc(1,0)
&
\Dt(b) =0= \cc(\nu,J)\\[5mm]
b\;=\; \young(11,24) \otimes \young(2,3)
&
\quad (\nu,J) \;=
& \yngrc(1,-1) & \yngrc(1,0,1,0) & \yngrc(1,0)
&
\Dt(b) =1= \cc(\nu,J)\\[5mm]
b\;=\; \young(12,23) \otimes \young(1,4)
&
\quad (\nu,J) \;=
& \yngrc(1,0) & \yngrc(1,0,1,0) & \yngrc(1,-1)
&
\Dt(b) =1= \cc(\nu,J)\\[5mm]
b\;=\; \young(12,24) \otimes \young(1,3)
&
\quad (\nu,J) \;=
& \yngrc(1,0) & \yngrc(1,0,1,-1) & \yngrc(1,0)
&
\Dt(b)=1= \cc(\nu,J)\\[5mm]
b\;=\; \young(13,24) \otimes \young(1,2)
&
\quad (\nu,J) \;=
& \yngrc(1,0) & \yngrc(1,0,1,0) & \yngrc(1,0)
&
\Dt(b)=2= \cc(\nu,J)\\[5mm]
b\;=\; \young(11,23) \otimes \young(2,4)
&
\quad (\nu,J) \;=
& \yngrc(1,-1) & \yngrc(2,0) & \yngrc(1,-1)
&
\Dt(b)=0= \cc(\nu,J)\\[5mm]
b\;=\; \young(12,34) \otimes \young(1,2)
&
\quad (\nu,J) \;=
& \yngrc(1,-1) & \yngrc(2,1) & \yngrc(1,-1)
&
\Dt(b)=1= \cc(\nu,J)
\end{array}
\end{equation*}
The unrestricted Kostka polynomial in this case is $M(L,\la;q)=2+4q+q^2=X(B,\la;q)$.
\end{example}

\subsection{Characterization of unrestricted rigged configurations}\label{sec:lower bound}
In this section we give an explicit description of the elements in
$\RC(L,\la)$ for type $A_{n-1}^{(1)}$. Generally speaking, the elements
are rigged configurations where the labels lie between the vacancy number
and certain lower bounds defined explicitly. This characterization will be used
in the next section to write down an explicit fermionic formula $M(L,\la;q)$
for the unrestricted configuration sum $X(B,\la;q)$.

Let $L=(L_i^{(a)} \mid (a,i)\in \HH)$ be a multiplicity array and 
$\la=(\la_1,\ldots,\la_n)$ be the $n$-tuple of nonnegative integers. 
The set of $(L,\la)$-configurations $\Conf(L,\la)$ is the set of all sequences 
of partitions $\nu=(\nu^{(a)} \mid a\in I)$ such that~\eqref{eq:constraint} holds.
As discussed in Section~\ref{sec:RC}, in the usual setting a rigged 
configuration $(\nu,J)\in \RCb(L,\la)$ consists of a configuration
$\nu\in \Confb(L,\la)$ together with a double sequence of partitions
$J=\{J^{(a,i)}\mid (a,i)\in\HH \}$ such that the partition
$J^{(a,i)}$ is contained in a $m_i^{(a)}\times p_i^{(a)}$ rectangle.
In particular this requires that $p_i^{(a)}\ge 0$. The unrestricted rigged
configurations $(\nu,J)\in \RC(L,\la)$ can contain labels that are negative,
that is, the lower bound on the parts in $J^{(a,i)}$ can be less than zero.

To define the lower bounds we need the following notation. 
Let $\la'=(c_1,c_2,\ldots,c_{n-1})^t$,
where $c_k=\la_{k+1}+\la_{k+2}+\cdots+\la_n$ is the length of the $k$-th
column of $\la'$, and let $\A(\la')$ be the set of tableaux of shape $\la'$ such that 
the entries are strictly decreasing along columns, and the letters in column $k$ are from 
the set $\{1,2,\ldots,c_{k-1}\}$ with $c_0=c_1$.

\begin{example} For $n=4$ and $\la=(0,1,1,1)$, the set $\A(\la')$
consists of the following tableaux
\begin{equation*}
\young(332,22,1) \quad \young(332,21,1) \quad \young(322,21,1) \quad 
\young(331,22,1) \quad \young(331,21,1) \quad \young(321,21,1).
\end{equation*}
\end{example}

\begin{remark}\label{rem:row}
Denote by $t_{j,k}$ the entry of $t\in\A(\la')$ in row $j$ and column $k$.
Note that $c_k-j+1\le t_{j,k}\le c_{k-1}-j+1$ since the entries in column $k$
are strictly decreasing and lie in the set $\{1,2,\ldots, c_{k-1}\}$.
This implies $t_{j,k}\le c_{k-1}-j+1\le t_{j,k-1}$, so that the rows of $t$
are weakly decreasing.
\end{remark}

Given $t\in\A(\la')$, we define the \textbf{lower bound} as
\begin{equation*}
M_i^{(a)}(t)=-\sum_{j=1}^{c_a} \chi(i\ge t_{j,a})
+\sum_{j=1}^{c_{a+1}} \chi(i\ge t_{j,a+1}),
\end{equation*}
where recall that $\chi(S)=1$ if the the statement $S$ is true and $\chi(S)=0$ otherwise.

Let $M,p,m\in \Z$ such that $m\ge 0$.
A $(M,p,m)$-quasipartition $\mu$ is a tuple of integers $\mu=(\mu_1,\mu_2,\ldots,\mu_m)$
such that $M\le \mu_m\le \mu_{m-1}\le \cdots\le \mu_1\le p$. Each $\mu_i$ is called
a part of $\mu$. Note that for $M=0$ this would be a partition with at most $m$ parts 
each not exceeding $p$.

The following theorem shows that the set of unrestricted rigged configurations
can be characterized via the lower bounds.
\begin{theorem} \cite[Theorem 4.6]{S:2005a}
Let $(\nu,J)\in \RC(L,\la)$. Then $\nu\in\Conf(L,\la)$ and
$J^{(a,i)}$ is a $(M_i^{(a)}(t),p_i^{(a)},m_i^{(a)})$-quasipartition
for some $t\in \A(\la')$. Conversely, every $(\nu,J)$ such that 
$\nu\in\Conf(L,\la)$ and $J^{(a,i)}$ is a 
$(M_i^{(a)}(t),p_i^{(a)},m_i^{(a)})$-quasipartition for some $t\in \A(\la')$
is in $\RC(L,\la)$.
\end{theorem}

\begin{example}
Let $n=4$, $\la=(2,2,1,1)$, $L_1^{(1)}=6$ and all other $L_i^{(a)}=0$. Then
\begin{equation*}
(\nu,J) \;=\; \yngrc(3,-2,1,0) \quad \yngrc(2,0) \quad \yngrc(1,-1)
\end{equation*}
is an unrestricted rigged configuration in $\RC(L,\la)$, where we have written
the parts of $J^{(a,i)}$ next to the parts of length $i$ in partition $\nu^{(a)}$.
To see that the riggings form quasipartitions, let us write
the vacancy numbers $p_i^{(a)}$ next to the parts of length $i$ in partition $\nu^{(a)}$:
\begin{equation*}
\yngrc(3,0,1,3) \quad \yngrc(2,0) \quad \yngrc(1,-1).
\end{equation*}
This shows that the labels are indeed all weakly below the vacancy numbers. For
\begin{equation*}
\young(441,33,2,1) \in \A(\la')
\end{equation*}
we get the lower bounds
\begin{equation*}
\yngrc(3,-2,1,-1) \quad \yngrc(2,0) \quad \yngrc(1,-1),
\end{equation*}
which are less or equal to the riggings in $(\nu,J)$.
\end{example}

For type $A_1$ we have $\la=(\la_1,\la_2)$ so that $\A=\{ t\}$
contains just the single tableau
\begin{equation*}
t=\begin{array}{|c|} \hline \la_2\\ \hline \la_2-1\\ \hline \vdots\\ \hline 1\\
\hline \end{array}.
\end{equation*}
In this case $M_i(t)=-\sum_{j=1}^{\la_2} \chi(i\ge t_{j,1})=-i$. This agrees
with the findings of~\cite{T:2005}.

As we will see in section~\ref{sec:level} the characterization of unrestricted
rigged configurations is similar to the characterization of level-restricted 
rigged configurations~\cite[Definition 5.5]{SS:2001}. Whereas the unrestricted
rigged configurations are characterized in terms of lower bounds, for level-restricted
rigged configurations the vacancy number has to be modified according to tableaux
in a certain set.

\subsection{Fermionic formula} \label{sec:fermi}
With the explicit characterization of the unrestricted rigged configurations
of Section~\ref{sec:lower bound}, it is possible to derive an explicit formula 
for the polynomials $M(L,\la)$ of~\eqref{eq:XM def}.

Let $\SA(\la')$ be the set of all nonempty subsets of $\A(\la')$ and set
\begin{equation*}
M_i^{(a)}(S)=\max\{M_i^{(a)}(t) \mid t\in S\} \qquad \text{for $S\in\SA(\la')$.}
\end{equation*}
By inclusion-exclusion the set of all allowed riggings for a given $\nu\in\Conf(L,\la)$ is
\begin{equation*}
\bigcup_{S\in\SA(\la')} (-1)^{|S|+1} \{J\mid \text{$J^{(a,i)}$ is a 
 $(M_i^{(a)}(S),p_i^{(a)},m_i^{(a)})$-quasipartition}\}.
\end{equation*}
The $q$-binomial coefficient $\qbin{m+p}{m}$, defined as
\begin{equation*}
\qbin{m+p}{m}=\frac{(q)_{m+p}}{(q)_m(q)_p},
\end{equation*}
where $(q)_n=(1-q)(1-q^2)\cdots(1-q^n)$, is the generating function of partitions
with at most $m$ parts each not exceeding $p$. Hence the polynomial $M(L,\la)$
may be rewritten as 
\begin{multline} \label{eq:fermi M}
M(L,\la;q)=\sum_{S\in\SA(\la')} (-1)^{|S|+1} \sum_{\nu\in\Conf(L,\la)}
q^{\cc(\nu)+\sum_{(a,i)\in\HH} m_i^{(a)}M_i^{(a)}(S)}\\
\times \prod_{(a,i)\in\HH} \qbin{m_i^{(a)}+p_i^{(a)}-M_i^{(a)}(S)}{m_i^{(a)}}
\end{multline}
called \textbf{fermionic formula}. By Corollary~\ref{cor:X=M} this is also a formula
for the unrestricted configuration sum $X(B,\la;q)$. This formula is different from the 
fermionic formulas of~\cite{HKKOTY:1999,Kir:2000}
which exist in the special case when $L$ is the multiplicity array of
$B=B^{1,s_k}\otimes \cdots \otimes B^{1,s_1}$ or $B=B^{r_k,1}\otimes \cdots 
\otimes B^{r_1,1}$.

\subsection{The Kashiwara operators $e_0$ and $f_0$}\label{sec:affine}
The Kirillov--Reshetikhin crystals $B^{r,s}$ are affine crystals and admit
the Kashiwara operators $e_0$ and $f_0$. As we have seen in~\eqref{eq:affine}
they can be defined in terms of the \textbf{promotion operator} $\pr$ as
\begin{equation*}
e_0=\pr^{-1}\circ e_1\circ \pr \quad \text{and} \quad 
f_0=\pr^{-1}\circ f_1\circ \pr.
\end{equation*}

We are now going to define the promotion operator on unrestricted rigged configurations.
\begin{definition}\label{def:rc pr}
Let $(\nu,J)\in\RC(L,\la)$. Then $\pr(\nu,J)$ is obtained as follows:
\begin{enumerate}
\item Set $(\nu',J')=f_1^{\la_1} f_2^{\la_2} \cdots f_n^{\la_n}(\nu,J)$ where
$f_n$ acts on $(\nu,J)^{(n)}=\emptyset$.
\item Apply the following algorithm $\rho$ to $(\nu',J')$ $\la_n$ times: Find the smallest
singular string in $(\nu',J')^{(n)}$. Let the length be $\ell^{(n)}$.
Repeatedly find the smallest singular string in $(\nu',J')^{(k)}$ of length
$\ell^{(k)}\ge \ell^{(k+1)}$ for all $1\le k<n$. Shorten the selected strings
by one and make them singular again.
\end{enumerate}
\end{definition}

\begin{example}
Let $B=B^{2,2}$, $L$ the corresponding multiplicity array and $\la=(1,0,1,2)$. Then
\begin{equation*}
(\nu,J)=\quad \yngrc(1,0) \quad \yngrc(2,-1,1,-1) \quad \yngrc(2,-1) \quad
\in \RC(L,\la)
\end{equation*}
corresponds to the tableau $b=\young(13,44)\in \Path(B,\la)$. After step (1)
of Definition~\ref{def:rc pr} we have
\begin{align*}
(\nu',J')&=\quad \yngrc(2,-1) \quad \yngrc(2,1,1,0) \quad \yngrc(2,-1,1,-1) \quad
\yngrc(2,-1).\\
\intertext{Then applying step (2) yields}
\pr(\nu,J)&=\quad \emptyset \quad \yngrc(1,0) \quad \yngrc(1,-1)
\end{align*}
which corresponds to the tableau $\pr(b)=\young(11,24)$.
\end{example}

\begin{lemma} \cite[Lemma 4.10]{S:2005a} \label{lem:pr}
The map $\pr$ of Definition~\ref{def:rc pr} is well-defined and 
satisfies~\eqref{eq:pr com} for $1\le a\le n-2$ and~\eqref{eq:weight rotation} 
for $0\le a\le n-1$.
\end{lemma}

Lemma 7 of~\cite{Sh:2002} states that for a single Kirillov--Reshetikhin crystal
$B=B^{r,s}$ the promotion operator $\pr$ is uniquely determined by~\eqref{eq:pr com} 
for $1\le a\le n-2$ and~\eqref{eq:weight rotation} for $0\le a\le n-1$. 
Hence by Lemma~\ref{lem:pr} $\pr$ on $\RC(L)$ is indeed the correct promotion
operator when $L$ is the multiplicity array of $B=B^{r,s}$.

\begin{theorem} \cite[Theorem 4.11]{S:2005a} \label{thm:pr}
Let $L$ be the multiplicity array of $B=B^{r,s}$. Then $\pr:\RC(L)\to\RC(L)$
of Definition~\ref{def:rc pr} is the promotion operator on rigged configurations.
\end{theorem}

\begin{conjecture} \cite[Conjecture 4.12]{S:2005a} \label{conj:pr}
Theorem~\ref{thm:pr} is true for any $B=B^{r_k,s_k}\otimes \cdots \otimes B^{r_1,s_1}$.
\end{conjecture}

Unfortunately, the characterization~\cite[Lemma 7]{Sh:2002} does not suffice to
define $\pr$ uniquely on tensor products $B=B^{r_k,s_k}\otimes \cdots \otimes B^{r_1,s_1}$.

\subsection{Open Problems} \label{sec:open X=M}
\begin{itemize}
\item In~\cite{DS:2005} a bijection $\Phi:\Path(B,\la)\to \RC(L,\la)$ is defined via a
direct algorithm. It is expected that Conjecture~\ref{conj:pr} can be proven by showing
that the following diagram commutes:
\begin{equation*}
\begin{CD}
\Path(B)@>{\Phi}>> \RC(L)\\
@V{\pr}VV @VV{\pr}V \\
\Path(B) @>>{\Phi}> \RC(L).
\end{CD}
\end{equation*}
Alternatively, an independent characterization of $\pr$ on tensor factors would give
a new, more conceptual way of defining the bijection $\Phi$ between paths and
(unrestricted) rigged configurations. A proof that the crystal operators
$f_a$ and $e_a$ commute with $\Phi$ for $a=1,2,\ldots,n-1$ is given in~\cite{DS:2005}.
\item Stembridge's local characterization of simply-laced crystals~\cite{St:2003} was used
in~\cite{S:2005a} to show that $f_a$ and $e_a$ of Definition~\ref{def:crystal} are 
in fact crystal operators. For nonsimply-laced types a local characterization of
crystals is not known yet. It can be shown via virtual crystals what the crystal
operators are in this case. See for example~\cite{OSS:2003a,OSS:2003b,S:2005b}
\item Hatayama et al.~\cite{HKKOTY:1999} derived a different fermionic formula
$M(L,\la;q)$ for the cases $B=B^{1,s_k}\otimes \cdots \otimes B^{1,s_1}$ and
$B=B^{r_k,1}\otimes \cdots \otimes B^{r_1,1}$. In~\cite{S:2002} this formula was interpreted
in terms of ``ribbon'' rigged configurations. It would be very interesting to 
relate the two fermionic formulas, in particular the two different rigged configurations.
As the fermionic formula of~\cite{HKKOTY:1999} is a special case of the Lascoux--Leclerc--Thibon
(LLT) spin generating function~\cite{LLT:1997}, this would yield a proof of a conjecture
by Kirillov and Shimozono~\cite[Conjecture 5]{KS:2002} that the LLT spin generating 
function labeled by a partition whose $k$-quotient is a sequence of rectangles is the 
same as the unrestricted generalized Kostka polynomial $X(B,\la;q)$.
\item The unrestricted rigged configurations for the $A_1$ case also appeared in a paper
by Takagi~\cite{T:2005} in the study of box-ball systems. A similar link should be
given for the general $A_{n-1}$ case.
\item Bailey's lemma is a powerful tool to prove Rogers--Ramanujan-type identities.
Andrews~\cite{Andrews:1984} showed that Bailey's lemma has an iterative structure
which relies on a transformation property of the $q$-binomial coefficients. This
iterative structure allows to derive infinite families of Rogers--Ramanujan identities
from a single seed identity. Since the unrestricted configuration sums $X$ yield
a generalization of the $q$-binomial coefficients, it is expected that they also
satisfy certain transformation properties which would give rise to a Bailey lemma.
For type $A_2$ this has been achieved in~\cite{ASW:1999}. The explicit formula $M$
for the unrestricted configuration sum might trigger further progress on
generalizations of the Bailey lemma to higher rank and other types.
\item For type $D_n^{(1)}$, a simple characterization in terms of lower bounds
for the parts of a configuration $\nu\in C(L)$ does not seem to exist.
For example take $B=B^{2,1}$ of type $D_4^{(1)}$ so that $L_1^{(2)}=1$ and all 
other $L_i^{(a)}=0$. Then the unrestricted rigged configurations
\begin{equation*}
\yngrc(1,0) \quad \yngrc(1,0,1,0) \quad \yngrc(1,0) \quad \yngrc(1,0)
\quad \text{and}\quad
\yngrc(1,0) \quad \yngrc(1,0,1,-1) \quad \yngrc(1,0) \quad \yngrc(1,0),
\end{equation*}
which correspond to the crystal elements $\begin{array}{|c|} \hline 1\\ \hline
\overline{1} \\ \hline \end{array}$ and $\begin{array}{|c|} \hline 3\\ \hline
\overline{3}\\ \hline \end{array}$ respectively, occur in $\RC(L)$, but
\begin{equation*}
\yngrc(1,0) \quad \yngrc(1,-1,1,-1) \quad \yngrc(1,0) \quad \yngrc(1,0)
\end{equation*}
on the other hand does not appear. It remains to determine a closed form
fermionic expression in this case.
\end{itemize}

\section{$X^\ell=M^\ell$}
\label{sec:level}

The fermionic formula for the level-restricted $X^\ell=M^\ell$ theorem has
a similar structure to the unrestricted fermionic formula. Instead of 
modifying the lower bounds for the rigged configurations, the upper bounds
are adapted.

\subsection{Level-restricted rigged configurations}
\label{sec:level rc}
A partition $\la=(\la_1,\la_2,\ldots,\la_n)$ is restricted of level $\ell$
if $\la_1-\la_n\le \ell$. Here $\la$ has at most $n$ parts, some of which may be zero.
Fix a shape $\la$ that is restricted of level $\ell$ and let $L$ be a multiplicity
array such that $L_i^{(a)}=0$ if $i>\ell$. Call such a multiplicity array
level-$\ell$ restricted. Define
$\lt = \ell-(\la_1-\la_n)$, which is nonnegative by assumption.

Set $\la'=(\la_1-\la_n,\ldots,\la_{n-1}-\la_n)^t$ and 
denote the set of all column-strict tableaux of shape $\la'$ over the alphabet
$\{1,2,\ldots,\la_1-\la_n\}$ by $\CST(\la')$. Define a table of modified
vacancy numbers depending on $\nu\in\Conf(L,\la)$ and $t\in\CST(\la')$ by
\begin{equation} \label{t vacancy}
  p_i^{(k)}(t) =
  p_i^{(k)} - \sum_{j=1}^{\la_k-\la_n} \chi(i\ge\lt+t_{j,k})
+ \sum_{j=1}^{\la_{k+1}-\la_n} \chi(i\ge\lt+t_{j,k+1})
\end{equation}
for all $i,k\ge1$ and $t_{j,k}$ is the $(j,k)$-th entry of $t$.
Finally let $\x_i^{(k)}$ be the largest part of the partition
$J^{(k,i)}$; if $J^{(k,i)}$ is empty set $x_i^{(k)}=0$.

\begin{definition}\label{def levrc}
Say that $(\nu,J)\in\RCb(L,\la)$ is restricted of level $\ell$
provided that
\begin{enumerate}
\item $\nu_1^{(k)} \le \ell$ for all $k$.
\item There exists a tableau $t\in\CST(\la')$, such that for every 
$i,k\ge 1$,
\begin{equation*}
  \x_i^{(k)} \le p_i^{(k)}(t).
\end{equation*}
\end{enumerate}
Let $\Conf^\ell(L,\la)$ be the set of all $\nu\in \Confb(L,\la)$ such that
the first condition holds, and denote by $\RC^\ell(L,\la)$ the set of 
$(\nu,J)\in\RCb(L,\la)$ that are restricted of level $\ell$. 
\end{definition}
Note in particular that the second condition requires that
$p_i^{(k)}(t)\ge 0$ for all $i,k\ge 1$.

\begin{example}\label{ex kir}
Let us consider Definition \ref{def levrc} for two classes of shapes $\la$
more closely:
\begin{enumerate}
\item \label{e rect}
Vacuum case: Let $\la=(a^n)$ be rectangular with $n$ rows.
Then $\la'=\emptyset$ and $p_i^{(k)}(\emptyset)=p_i^{(k)}$
for all $i,k\ge 1$ so that the modified vacancy numbers are equal
to the vacancy numbers.
\item Two-corner case: Let $\la=(a^\alpha,b^\beta)$ with $\alpha+\beta=n$
and $a>b$. Then $\la'=(\alpha^{a-b})$ and there is only one tableau $t$
in $\CST(\la')$, namely the Yamanouchi tableau of shape $\la'$.
Since $t_{j,k}=j$ for $1\le k\le \alpha$ we find that
\begin{equation*}
p_i^{(k)}(t)=p_i^{(k)}-\delta_{k,\alpha}\max\{i-\lt,0\}
\end{equation*}
for $1\le i\le \ell$ and $1\le k<n$.
\end{enumerate}
\end{example}

We define the level-restricted rigged configuration generating function as
\begin{equation}\label{eq:level rc}
  M^\ell(L,\la;q) = \sum_{(\nu,J)\in\RC^\ell(L,\la)} q^{\cc(\nu,J)}.
\end{equation}
The $X^\ell=M^\ell$ conjecture was proven in~\cite{SS:2001}.
\begin{theorem} \cite[Theorem 5.7]{SS:2001}
For a level-$\ell$ restricted partition $\la$ and a level-$\ell$ restricted multiplicity
array $L$ we have $X^\ell(B,\la;q)=M^\ell(L,\la;q)$.
\end{theorem}

\begin{example} Consider $n=3$, $\ell=2$, $\la=(3,2,1)$, $L_1^{(1)}=4$, $L_2^{(1)}=1$
and all other $L_i^{(a)}=0$. Then
\begin{equation}\label{conf ex}
 \yngrc(1,0,1,0,1,0) \qquad \yngrc(1,1)
\qquad \text{and} \qquad
 \yngrc(2,1,1,2) \qquad \yngrc(1,0)
\end{equation}
are in $\Conf^\ell(L,\la)$ where again the vacancy numbers are indicated 
to the left of each part.
The set $\CST(\la')$ consists of the two elements
\begin{equation*}
 \young(11,2) \qquad \text \qquad \young(12,2).
\end{equation*}
Since $\lt=0$ the three rigged configurations
\begin{equation*}
 \yngrc(1,0,1,0,1,0) \quad \yngrc(1,0), \qquad
 \yngrc(2,0,1,0) \quad \yngrc(1,0) \qquad \text{and} \qquad
 \yngrc(2,0,1,1) \quad \yngrc(1,0)
\end{equation*}
are restricted of level 2 with charges $2,3,4$, respectively.
The riggings are given on the right of each part. 
Hence $M^\ell(L,\la;q)=q^2+q^3+q^4$.

In contrast to this, the rigged configuration generating function
$\Mb(L,\la;q)$ is obtained by summing over both configurations in~\eqref{conf ex}
with all possible riggings below the vacancy numbers.
This amounts to $\Mb(L,\la;q)=q^2+2q^3+2q^4+2q^5+q^6$.
\end{example}

\subsection{Level-restricted fermionic formula}
\label{sec:level fermi}
Similarly to the unrestricted case of section~\ref{sec:fermi}, one can rewrite the expression
of the level-restricted rigged configuration generating function of~\eqref{eq:level rc}
in fermionic form. It was shown in~\cite[Lemma 6.1]{SS:2001} that 
$p_i^{(k)}(t)=0$ for all $i\ge \ell$.

Let $\SCST(\la')$ be the set of all nonempty subsets of $\CST(\la')$.
Furthermore set $p_i^{(k)}(S)=\min\{p_i^{(k)}(t)|t\in S\}$ for
$S\in\SCST(\la')$. Then by inclusion-exclusion the set of allowed
rigging for a given configuration $\nu\in\Conf^\ell(L,\la)$ is given by
\begin{equation*}
\sum_{S\in\SCST(\la')}(-1)^{|S|+1} \{J|\x_i^{(k)}\le p_i^{(k)}(S)\}.
\end{equation*}
Since the $q$-binomial $\qbins{m+p}{m}$ is the generating function of 
partitions with at most $m$ parts each not exceeding $p$ and since
$p_\ell^{(k)}(S)=0$ by~\cite[Lemma 6.1]{SS:2001} the level-$\ell$ restricted 
fermionic formula has the following form.
\begin{theorem} \cite[Theorem 6.2]{SS:2001}
\begin{equation*}
M^\ell(L,\la;q)=\sum_{S\in\SCST(\la')} (-1)^{|S|+1}
\sum_{\nu\in\Conf^\ell(L,\la)} q^{\cc(\nu)}
\prod_{i=1}^{\ell-1}\prod_{k=1}^{n-1} 
\qbin{m_i^{(k)}+p_i^{(k)}(S)}{m_i^{(k)}}.
\end{equation*}
\end{theorem}

\subsection{Open Problems}

\begin{itemize}
\item
In~\cite[Conjecture 8.3]{SS:2001} it was conjectured that the bijection
$\Phib$ is also well-behaved with respect to fixing certain subtableaux
in the set of Littlewood-Richardson tableaux. In the crystal language
let $\rho\subset \la$ be a partition and 
$b_\rho=b_k\otimes \cdots \otimes b_1 \in B_\rho=B^{\rho_k^t,1} \otimes 
\cdots \otimes B^{\rho_1^t,1}$ where $b_i$ is the column tableau of height $\rho_i^t$ 
with $\row(b_i)=\rho_i^t \ldots 21$.
Denote the set of all paths in $\Path^\ell(B\otimes B_\rho,\la)$ with fixed subpath
$b_\rho$ by $\Path^\ell(B,\la,\rho)$.
Set $\rho'=(\rho_1-\rho_n,\ldots,\rho_{n-1}-\rho_n)^t$ and
\begin{equation*}
M_i^{(k)}(t)=\sum_{j=1}^{\rho_k-\rho_n}\chi(i\le \rho_1-\rho_n-t_{j,k})
-\sum_{j=1}^{\rho_{k+1}-\rho_n}\chi(i\le \rho_1-\rho_n-t_{j,k+1})
\end{equation*}
for all $t\in\CST(\rho')$.
Then define $\RC^\ell(L,\la,\rho)$ to be the set of all
$(\nu,J)\in\RC^\ell(L\cup L_\rho,\la)$ such that there exists
a $t\in\CST(\rho')$ such that $M_i^{(k)}(t)\le x$
for $(i,x)\in (\nu,J)^{(k)}$ and $M_i^{(k)}(t)\le p_i^{(k)}$
for all $i,k\ge 1$. Here $L_\rho$ is the multiplicity array of $B_\rho$.
Note that the second condition is obsolete
if $i$ occurs as a part in $\nu^{(k)}$ since by definition
$M_i^{(k)}(t)\le x\le p_i^{(k)}$ for all $(i,x)\in (\nu,J)^{(k)}$.
Conjecture 8.3 of~\cite{SS:2001} asserts that $\Path^\ell(B,\la,\rho)$ and 
$\RC^\ell(L,\la,\rho)$ correspond under $\Phib$.
\item It is still an open problem to provide a combinatorial formula
for the fusion coefficients of the Verlinde algebra~\cite{TUY:1989,V:1988}.
The fermionic formulas of this section only provide such a formula for
rectangular tensor factors.
\end{itemize}

\end{document}